\documentclass[12pt]{amsart}
\usepackage{amscd,graphicx}
\newcommand{\bi}{{\bf{i}}}
\newcommand{\bj}{{\bf{j}}}
\newcommand{\bk}{{\bf{k}}}
\newcommand{\vP}{\vec{P}}

\newcommand{\g}{su(2)}
\newcommand{\G}{SU(2)}
\newcommand{\vv}{\vec{v}}
\newcommand{\vw}{\vec{w}}
\newcommand{\vu}{\vec{u}}
\newcommand{\bZ}{\mathbb{Z}}
\newcommand{\bR}{\mathbb{R}}
\newcommand{\bC}{\mathbb{C}}

\newcommand{\kg}{\pi_1(S^3-K)}

\newtheorem{theorem}{Theorem}
\newtheorem{lemma}{Lemma}
\newtheorem{prop}{Proposition}

\newtheorem{defi}{Definition}
\newtheorem{conj}{Conjecture}

\title{Dubois' Torsion, A-Polynomial and Quantum Invariants.}
\author{Charles Frohman}

\address{Department of Mathematics, The University of Iowa, Iowa City, IA
52242, USA}

\email{\tt charles-frohman@.uiowa.edu}

\author{Joanna Kania-Bartoszynska}

\address{National Science Foundation, Arlington, VA, 22230, USA} 

\email{\tt jkaniaba@nsf.gov}
\thanks{This material is based upon work supported by and while serving at the National Science Foundation. Any opinion, findings, and conclusions or recommendations expressed in this material are those of the authors and do not necessarily reflect the views of the National Science Foundation.}

\begin{document}

\begin{abstract}It is shown that for knots with a sufficiently regular character variety the Dubois' torsion detects the $A$-polynomial of the knot. A global formula for the integral of the Dubois torsion is given. The formula
looks like the heat kernel regularization of the formula for the Witten-Reshetikhin-Turaev invariant of the double of the knot complement.   The Dubois' torsion is recognized
as the pushforward of a measure on the character variety of the double of the knot complement coming from the square root of Reidemeister torsion.  This is used to motivate a conjecture about quantum invariants detecting the $A$-polynomial. \end{abstract}

\maketitle

\section{Introduction} There is a close connection between the colored Jones polynomial and
the $SL_2{\mathbb{C}}$ characters of knot groups.  The first intimation of this connection
appeared in Bullock's work on the skein modules of knot complements \cite{B1}.  The idea was developed more formally in \cite{B} and \cite{PS}. The point is that the Kauffman bracket skein module of a three manifold $M$, when the complex parameter is set to be equal to $-1$, is the coordinate ring of the unreduced scheme of the
$SL_2\mathbb{C}$-characters of the fundamental group of $M$.
The connection with the $A$-polynomial of Cooper, Culler, Gillet, Long and Shalen \cite{CCLS} was first developed in \cite{FG} and
\cite{FGL}. The idea was to first realize the $A$-polynomial as generating the kernel of a map between rings of characters, and then deform the rings in a canonical way so that the kernel can be understood
as a submodule of a skein module, called the $A$-ideal.   Garoufalidis and Le  developed a more algebraic setting for these ideas, that allowed them to prove nontriviality of the non-commutative $A$-ideal, along with making many compelling conjectures \cite{LS}.

It has long been conjectured that the asymptotics of the colored Jones polynomials of a knot
are modulated by the $SL_2{\mathbb{C}}$ or $SU(2)$-characters of the knot group. As $SL_2{\mathbb{C}}$ is the complexification of $SU(2)$ there is a common ground for trying to understand the relations between the representation theory of knot groups and the values
of the colored Jones polynomial.

 There is a  conjectured  
asymptotic formula for the Witten-Reshetikhin-Turaev invariants for links in $3$-manifolds \cite{W1,FrG,Je,JW,JW1,Roz}.  It is presented here as it appeared in \cite{Oh}, where it was  stated  only for the empty link in a closed oriented three manifold $M$.
\begin{equation}\label{wittenexp}
Z_r^{SU(2)}(M)\sim_{r\rightarrow \infty} e^{-3\pi \bf{i}(1+b^1(M))/4} 
\sum_{[A]}e^{2\pi \bf{i}CS(A)} r^{(h^1_A-h^0_A)/2} e^{-2\pi\bf{i}(I_A/4+h^0_A/8)}\tau_M(A)^{1/2}
\end{equation}
Here $b^1(M)$ is the first Betti number of $M$, $CS(A)$ is the Chern-Simons invariant of the flat
connection $A$, and $[A]$ denotes the gauge equivalence class of $A$. Next, $h^i_A$ is the rank of the $i$-th cohomology of $M$ with coefficients in the $su(2)$-bundle twisted by the adjoint action of the monodromy of $A$, and $I_A$ is the spectral flow of the 
signature operator along a path connecting $A$ to the trivial flat connection.  Finally,
$\tau_M(A)$ is the Reidemeister torsion of the chain complex $C^i(M,ad A)$.

When the $SU(2)$-character variety has positive dimension, the sum in the above formula
is an integral.  In this case, if we restrict to connections that are irreducible,  the square root of the Reidemeister torsion defines a natural measure on the character variety. In the case of torus knots, Kashaev and Dubois have worked out the asymptotics of the colored Jones polynomial and found in \cite{DK} an asymptotic expansion having the same flavor as the formula (\ref{wittenexp}).

In this paper we explore the Reidemeister torsion of Dubois \cite{Du} as a measure on the regular part of the $SU(2)$-character variety of a knot complement.  Exploiting the connection with the Kauffman bracket skein algebra of the torus \cite{FG}, we use it to define a seminorm on that algebra. When the $SU(2)$-character variety is sufficiently nonsingular the
radical of this seminorm is the ideal of functions that vanish on the image of the irreducible
representations of the knot complement in the character variety of the torus.

We go on to develop a global formula for the seminorm that looks like the Witten-Reshetikhin-Turaev invariant of the skein in the double of the knot complement. This leads us to
a conjectural characterization of the $A$-polynomial of knots with sufficiently regular character varieties in terms of quantum invariants.

 Rozansky \cite{Roz} first recognized that for  regular representations of $3$-manifolds, the square root of the Reidemeister torsion defines an invariant volume form on the regular part of the character variety.
As part of the motivation of our conjecture  we prove that the Dubois' torsion of a knot is a geometrically motivated evaluation of the square root of the Reidemeister torsion of the double of the knot complement. 

The plan of the paper is as follows.  In section 2 we talk about the representation spaces of knot groups. In section 3 we review Dubois' torsion and describe the corresponding seminorm. In section 4 we give a global computation of the seminorm. In section 5 we compare the answer from the global computation to the formula for the Witten-Reshetikhin-Turaev invariants of links in the double
of the knot complement. In the last section we prove that DuboisÕ torsion of a knot complement is a geometrically motivated evaluation of the square root of the Reidemeister torsion of its double.
Throughout the paper we work the example of the trefoil knot to make the exposition concrete.

The authors thank Thomas Kerler who helped us understand the global computation, and Paul Kirk who helped with the interpretation of the conjectured asymptotic formula.

\section{Representations of knot groups}
In this section we recall the definition of the A-ideal of a knot and illustrate it with the computation for the trefoil.
\subsection{The $A$-ideal}\label{Aideal}
Let $K\subset S^3$ be a knot, and $N(K)$ an open regular neighborhood of $K$ so that $S^3-N(K)$ is a smooth manifold
with boundary a torus $T^2$. The fundamental group of the torus is $\bZ\times \bZ$ with generators
the longitude $\lambda$ and meridian $\mu$.  We are interested in the  image of the representations of the fundamental
group of $S^3-K$ in the representations of the fundamental group of  $T^2$  under restriction. Although, the study of images like this was initiated using $SL(2,\bC)$ representations \cite{CCLS}, in this paper we study  $\G$ representations.  In \cite{CCLS} the authors also pass to a  cover, so that the image is cut out by a single polynomial, which they dubbed the $A$-polynomial. It is common to throw out the component of the abelian representations
or just focus on the component that contains the holonomy of the complete hyperbolic structure and still refer
to the polynomial cutting out that variety as the $A$-polynomial.

Recall that $\G$ consists of two by two complex matrices of the form \begin{equation}\begin{pmatrix} \alpha & \beta \\ -\overline{\beta} &\overline{\alpha}\end{pmatrix}\end{equation} where $|\alpha|^2+|\beta|^2=1$. Alternatively, you can think of $SU(2)$ as the unit quaternions. In the quaternionic model, each element of $SU(2)$ can be written as $\cos{\phi}+\sin{\phi}\vP$, where $\vP$ is a unit vector
in $\bR^3$. This second model yields an elegant way of working with the tangent space of $SU(2)$ at the identity, that is,  \begin{equation}\g=T_{\mathbf{1}}S^3= \bR^3.\end{equation} Under this identification, the adjoint representation $Ad:SU(2)\rightarrow End(\g)$,
\begin{equation} Ad_X(\vv)=X\vv X^{-1}\end{equation} is easily
understood.   If $X=\cos{\phi}+\sin{\phi}\vP$ then $Ad_X$ acts as a rotation by $2\phi$ radians about the axis $\vP$.  This simple
picture led to the beautiful computations of $SU(2)$ representations of knot groups in \cite{Bu,K}.  To simplify
notation we denote the adjoint action of $X$ on $su(2)$ by a lower dot, 
$Ad_X(\vv)=X.\vv$.  The Lie bracket in $su(2)$ is twice the cross product. Given a representation of a group $\Gamma$, $\rho:\Gamma\rightarrow \G$,  denote by $Ad\rho$ the representation into
endomorphisms of $\g$ obtained by following $\rho$ with the adjoint action.

 If $M$ is a manifold with finitely generated fundamental group,
denote the space of representations of $\pi_1(M)$ into $SU(2)$ by $R(M)$.  The space $R(M)$  is realized as a subset of the Cartesian product $SU(2)^k$ where the coordinates are the values of the
representation on a finite set of generators of $\pi_1(M)$.  Those tuples that satisfy the equations coming from the group 
relators exactly correspond to $R(M)$, so that not only is $R(M)$ a topological space with the subspace topology
coming from $SU(2)^k$, it is also a real algebraic variety.  Actually, more is true:  by interpreting the group relators as matrix
equations, and then seeing each coefficient of the matrices as a real equation, the ideal generated by requiring the value
of each relator to be the identity is an invariant of $M$, see \cite{LM}.  The radical of that ideal is more commonly studied as an invariant. 

The group $SU(2)$ acts 
on $R(M)$ by conjugation. We denote the quotient space by $X(M)$. It inherits
the quotient topology from $R(M)$, and is also a real algebraic variety called the
{\em character variety}. The coordinate rings of 
these varieties are denoted by $C[R(M)]$ and $C[X(M)]$.  The ring $C[X(M)]$ can be identified with the subring $C[R(M)]^{\G}$ of $C[R(M)]$ that is fixed under the action of $SU(2)$. The ring $C[X(M)]$ is closely related to the Kauffman Bracket
skein module of $M$, as shown in \cite{B,BFK,PS}.

Any representation $\rho:\pi_1(T^2)\rightarrow \G$  is determined by the values on the longitude and meridian, $\rho(\lambda)$ and
$\rho(\mu)$ .  Since $\lambda$  and $\mu$ commute, the matrices $\rho(\lambda)$ and $\rho(\mu)$ are simultaneously
diagonalizable. Thus they are conjugate to a pair
\begin{equation}\label{deflandm}
(\begin{pmatrix} l & 0 \\ 0 & l^{-1}\end{pmatrix},\begin{pmatrix} m & 0 \\ 0 & m^{-1} \end{pmatrix}),\end{equation} 
where $l,m \in S^1\subset \bC$.  These coordinates are  ambiguous, since the two matrices
can be simultaneously conjugated by $\begin{pmatrix} 0 & -1 \\ 1 & 0\end{pmatrix}$ with the result of exchanging the positions of $l$ and $l^{-1}$ and $m$ and $m^{-1}$.  From this we see that the conjugacy classes of $\G$-representations of $\pi_1(T^2)$ can be identified with 
$S^1\times S^1/\sim$ ,where the equivalence relation $\sim$ comes from the hyperelliptic
involution $\theta:S^1\times S^1 \rightarrow S^1\times S^1$,  given by $\theta(l,m)=(l^{-1},m^{-1})$.  We have characterized $X(T^2)$ as a space that is called the {\em pillowcase} since it can be represented as the result of identifying two squares along their boundaries.

\vspace{.1in}

\begin{center}
\includegraphics{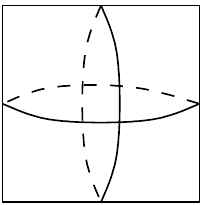}
\end{center}

\vspace{.1in}

 The coordinate ring of the pillowcase, $C[X(T^2)]$, can be understood
as the subring of ${\mathbb{C}}[l,l^{-1},m,m^{-1}]$ that is fixed by the endomorphism $\theta$.  It
can be shown that the functions of the form $l^pm^q+l^{-p}m^{-q}$ form a basis for $C[X(T^2)]$,  see \cite{FG}.

The inclusion map $i:T^2\rightarrow S^3-K$ induces restriction maps $R(S^3-K)\rightarrow R(T^2)$ and
$\psi:X(S^3-K)\rightarrow X(T^2)$.  We are especially interested in the ring homomorphism  induced by $\psi$, 
\begin{equation} \Psi :C[X(T^2)]\rightarrow C[X(S^3-K)].\end{equation}  By definition, if $f\in C[X(T^2)]$ then $\Psi(f)=f\circ \psi$.
 The radical of the kernel of $\Psi$ is an ideal $B(K)\subset C[X(T^2)]$
which we call the {\em $B$-ideal} of the knot $K$.  It is the ideal of polynomial functions on $X(T^2)$ that are zero
on the image of the restriction map $\psi$.  We are less interested in the characters of the abelian representations.
Let $A(K)$ be the ideal of functions that vanish on the image of the irreducible characters.  The study of the noncommutative
analog of this ideal was initiated in \cite{FGL}. The more sophisticated algebra of \cite{LS}, led to a simpler ideal
that was principle and could thus be described by a single polynomial. The nontriviality of that ideal is proved by
a beautiful argument based on the algebra of $q$-hypergeometric functions.

\subsection{An example}
The complement
of the trefoil collapses onto a $CW$-complex with one vertex $v$, two edges,
$x$ and $y$ and a single two cell whose attaching map is $x^2y^{-3}$.

\hspace{2in}\includegraphics{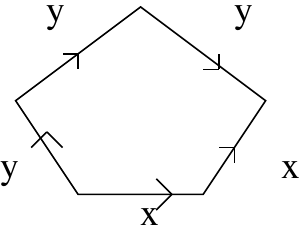}

This can be visualized by seeing the trefoil as lying on an unknotted
torus. In this setting there is a two complex consisting of the two cores of the
handlebodies bounded by the torus, and a singular annulus that runs
twice around one core and three times around the other, and misses the trefoil.  Pinch the singular annulus along a crosscut. The pinch point is the vertex, the remnants of the cores are the two edges, and the result of pinching the annulus is the two-cell.

This two--complex gives rise to a presentation of the fundamental group of the knot complement. There is one generator for each edge and one relator for each two-cell.  It is important to be able to express the meridian and longitude in terms of the generators.  The complement of the trefoil knot
has the fundamental group 
\begin{equation}\pi_1(S^3-K)=<x,y|x^2=y^3>.\end{equation}
The meridian is given by $\mu=xy^{-1}$, and the longitude is $\lambda=x^2\mu^{-6}=x^2(xy^{-1})^{-6}$.
The merdian and longitude commute, hence the word $\lambda\mu\lambda^{-1}\mu^{-1}$ in the free group on $x$ and $y$ is in the normal closure of the relator $r=x^{-2}y^3$.
Following Dubois \cite{DuF},
\begin{equation}\lambda\mu\lambda^{-1}\mu^{-1}=x^2\mu^{-6}\mu\mu^6x^{-2}\mu^{-1}=xrx^{-1}\  \mu r^{-1}\mu^{-1}.\end{equation}  From this equation, we can build a map of the torus $T^2$ into the two-complex of the knot complement  which is homotopic to the inclusion of the boundary torus.  Fill in a square with sides from the commutator by two $2$-cells corresponding to the word above. This gives a map of a disk into the two-complex, so that points on opposite sides of the square get mapped to the same point. The map descends to a map of a torus into the two-complex. It is clear that this torus represents the boundary of the knot complement.

The representations of $\pi_1(S^3-K)$ into $\G$  can be identified with a subset of  $SU(2)\times SU(2)$. Give $SU(2)\times SU(2)$ coordinates $(X,Y)$ and send $\rho:\kg \rightarrow \G$ to $(\rho(x),\rho(y))$.  To correspond to a representation, a pair $(X,Y)$ must satisfy the equation $X^2=Y^3$.  There are two components of this subset of $SU(2)\times SU(2)$: the points coming from
abelian representations, and the points coming from irreducible representations.

From \cite{K} we know that if
$\rho:\pi_1(S^3-K)\rightarrow \G$ is irreducible then $\rho(x^2)=\rho(y^3)=-1$.  This leads to an easy parametrizations of the components of $R(S^3-K)$.
To parametrize the abelian representations let $\alpha:SU(2)\rightarrow SU(2)\times SU(2)$ be given by $\alpha(Z)=(Z^3,Z^2)$.
To parametrize the irreducible representations, let
$\beta:SU(2)\times (0,\pi) \rightarrow SU(2)\times SU(2)$
be given by 
\begin{equation}\label{beta} 
\beta (A,t)=(A{\bf{i}}A^{-1},A(\cos{\frac{\pi}{3}}+\sin{\frac{\pi}{3}}(\cos{t}{\bf{i}}+\sin{t}{\bf{j}}))A^{-1}).\end{equation} 
The parametrization $\alpha$ of the abelian representations is one to one. The parametrization $\beta$  of the irreducible representations is $2$ to $1$, because the center of $SU(2)$ is in the kernel of the map coming from conjugation.
Notice that the abelian representations have codimension $3$ in $SU(2)\times SU(2)$ and the
irreducible representations have codimension $2$. 

Notice that the map $\beta$ extends smoothly to 
\begin{equation}\overline{\beta}:SU(2)\times [0,\pi]\rightarrow SU(2)\times SU(2).\end{equation}
The images of $SU(2)\times \{0\}$ and $SU(2)\times \{\pi\}$ under $\overline{\beta}$ are two spheres of abelian
representations.  The points in the abelian representations that lie in the closure of the irreducible representations are called {\em bifurcation points}. At the bifurcation points the squares of the eigenvalues of the matrices that are images of the meridian are roots of the Alexander polynomial of the trefoil knot. This phenomenon occurs for any knot \cite{D,K,FK}.  

The image of the irreducible representations  of the complement of the trefoil in the pillowcase under restriction  is depicted in the figure below.

\vspace{.1in}

\begin{center}
\includegraphics{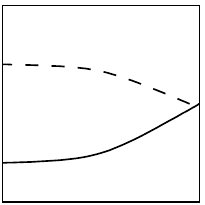}
\end{center}

\vspace{.1in} 
The $A$-polynomial of the trefoil is given by
\begin{equation}
A(l,m)= l+m^{-6}
\end{equation}
where $l$ and $m$ are the upper left hand corner of the diagonalized $\rho(\lambda)$ and $\rho(\mu)$ respectively, as in  (\ref{deflandm}).
The $A$-ideal of the trefoil is generated in ${\mathbb{C}}[l,l^{-1}, m,m^{-1}]$ by the polynomials 
\begin{equation}
l+l^{-1}+m^6+m^{-6}
\end{equation}
and
\begin{equation}
lm+l^{-1}m^{-1}+m^5+m^{-5}
\end{equation}
The computation can be found in \cite{Gelca}.

\section{Dubois' torsion}

 This section presents Dubois torsion for knots, which in turn is used  to construct a seminorm on $C[X(T^2)]$ whose radical is $A(K)$.

\subsection{Reidemeister torsion and volume}

If ${\bf{u}}=\{\vu_i\}$ and ${\bf{w}}=\{\vw_j\}$ are two ordered bases for the finite dimensional vector space $V$, then
each $\vw_j=\alpha^i_j\vu_i$, where, following Einstein's summation convention, the last equation represents a sum over all $i$.  The determinant of 
the change of basis matrix $(\alpha^i_j)$ is denoted $[{\bf{w}/\bf{u}}]$. 
Given a finite dimensional chain complex $C=(C^i,\delta^{i+1})$ whose chain
groups are vector spaces over $\bR$ with preferred ordered basis ${\bf{c}}_i$ for $C^i$, and 
a collection ${\bf{h}}_i$ of $i$ chains  that give rise to a basis for $H^i(C)$, choose an ordered basis 
	${\bf{b}}_i $ for the complement of $Z^i(C)$. Finally let ${\bf{\overline{b}}}_i$ be the image of ${\bf{b}}_i$ under $\delta^{i+1}:C^i \rightarrow C^{i+1}$.
Notice that ${\bf{\overline{b}}}_{i-1},{\bf{h}}_i,{\bf{b}}_i$ is an ordered basis for $C^i$.  The {\em Reidemeister torsion} of the chain complex $C$ evaluated on the basis ${\bf h}$ is defined by:
\begin{equation}\label{reidtor}
\tau({\bf{h}})=\left|\frac{\prod_{i \ odd} [{\bf{\overline{b}}}_{i-1},{\bf{h}}_i,{\bf{b}}_i/{\bf{c_i}}]}{\prod_{i \ even}[{\bf{\overline{b}}}_{i-1},{\bf{h}}_i,{\bf{b}}_i/{\bf{c_i}}]}\right|.
\end{equation}
In the case where the chain groups are innerproduct spaces, there is a preferred class of bases. Specifically, any orthonormal basis can be used.
This is because the change of basis matrix between two orthonormal bases has determinant $\pm 1$. Since we are working with the absolute value,  $\tau(h)$ does not depend on the choice of bases for $C^i$  as long as ${\bf{c}}_i$ are orthonormal.

The Reidemeister torsion transforms as a top dimensional form, i.e. a volume form,
on 
\begin{equation} \oplus_{i \ even} H^i(C) \bigoplus \oplus_{i \ odd} (H^i(C))^*.\end{equation}
The definition and properties of Reidemeister torsion can be found in \cite{M,RS,T1,T2,Wi}.
The tools developed for manipulating torsion in \cite{M} are fundamental.

This paper relies on interpreting the Reidemeister torsion of a chain complex in terms of Riemannian volume forms. To make this interpretation clear, we review some geometry of volume. Given an oriented innerproduct space $V$ there is a distinguished top dimensional form characterized by the fact that it evaluates to $1$ on any positively oriented orthonormal ordered basis of $V$. We call this the {\em volume form}. As forms are about determinants, the fact that the basis was orthonormal was not particularly important. 

Recall that, if $L:V\rightarrow W$ is a linear map of innerproduct spaces, we can restrict the domain of $L$ to the perpendicular
of its kernel and restrict its range to its image, $L:(kerL)^{\perp}\rightarrow imL$, to get a map of innerproduct spaces. We don't have orientations for these two spaces, so we don't know the sign of the determinant of $L$, but we can define $||L||$ to
be 
\begin{equation}\label{norm}
||L||=\sqrt{det(L^*\circ L)}=\sqrt{det(L\circ L^*)},
\end{equation} 
where we have restricted $L$ as above.   

A choice of vectors $\vv_1,\ldots,\vv_n$ of an innerproduct space
$V$ is equivalent to defining a linear map $L:\mathbb{R}^n\rightarrow V$.  The norm
of $L$, 
$||L||$, can be thought of as the $n$-dimensional content of the parallelepiped spanned by the vectors $\vv_i$.
Given a single vector $\vv\in V$ let $L:\mathbb{R}\rightarrow V$ be an operator defined by 
$L(1)=\vv$. Notice that $||L||=||\vv||$. Thus even though there is potentially some ambiguity in the use of the double bar notation,  the value is the same regardless of the meaning.

Let $M$ be an oriented Riemannian $m$-manifold. The pointwise choice of volume forms defines a smooth $m$-form on $M$, $\nu^M$, called the {\em Riemannian volume form}.  If $N\subset M$ is a smooth oriented submanifold, then the Riemannian metric on $M$ induces a Riemannian metric on $N$ and  an associated volume form $\nu^N$.  
The volume forms $\nu^M$ and $\nu^N$ are related as follows. 
At any $p\in N$ choose an orthonormal basis $n_1,\ldots,n_k$ for the orthogonal complement of $T_pN$ in $T_pM$ so that
as oriented vector spaces,
\begin{equation}<n_1,\ldots, n_k> \oplus T_pN =T_pM.\end{equation}
The Riemannian volume form on $N$ is the interior product of $\nu^M$ with $\{n_1,\ldots,n_k\}$. That is, for any $(w_1,\dots w_n)\in (T_p N)^n$,

\begin{equation} \nu^N(w_1,\ldots,w_n)=\nu^M(n_1,\ldots,n_k,w_1,\ldots,w_n).\end{equation}

Suppose instead that $v_1,\ldots,v_k$ is an arbitrary basis (not necessarily orthonormal) for the orthogonal complement
of $T_pN$. Then
\begin{equation} \nu^M(v_1,\ldots,v_k,w_1,\ldots,w_n)=\pm ||L||\nu^N(w_1,\ldots,w_n),
\end{equation} 
where $L:\mathbb{R}^k\rightarrow T_pM$ is the linear map that takes the standard basis for $\mathbb{R}^k$ to the vectors $v_1,\ldots,v_k$.

\subsection{Local definition of Dubois torsion}
Let $K$ be a knot and let $W \subset S^3-N(K)$ be a two-complex, with one vertex, embedded in $S^3-K$ so
that $S^3-N(K)$ is a closed regular neighborhood of $W$. Alternatively, one can think of $S^3-N(K)$ as having a
handle structure with a single zero-handle, $k$ one-handles, and $k-1$ two-handles.  Denote the vertex by
$v_1$, the edges by $e_i$ and the $2$-cells by $f_j$.  There is a presentation of the fundamental group of $S^3-K$ corresponding to this complex, with one generator $x_i$ for each edge $e_i$, and one relator $r_j$ for each face $f_j$. The relator $r_j$ is a word
in the generators $x_i$ coming from the attaching map of the cell $f_j$.
Let $C_*(W)$ be the chain complex that corresponds to the homology of
$W$ with basis $v_1,e_i,f_j$. The dual chain complex $C^*(W)=Hom(C_*(Y),\mathbb{R})$
has the dual basis $v^1,e^i,f^j$. For instance $e^i:C_*(W)\rightarrow \mathbb{R}$
is the linear map so that $e^i(e_i)=1$ and the value of $e^i$ applied to any other basis element is $0$.
We make $C^*(W)$ into an innerproduct space by declaring the basis $v^1,e^i,f^j$ to 
be orthonormal.

Define a map
 \begin{equation}\label{complex} 
 \overline{r}:SU(2)^k\rightarrow SU(2)^{k-1}
 \end{equation}
 by instantiating the variables $x_i$ in the relators $(r_1,\ldots, r_{k-1})$ by the matrices
$(X_1,\ldots,X_k)$. The representation variety $R(S^3-K)$ can be identified
with $\overline{r}^{-1}(Id,\ldots,Id)$.  The group $\G$ acts on $SU(2)^k$ by simultaneous conjugation of the entries, $c: SU(2)\times SU(2)^k\rightarrow SU(2)^k$,
\begin{equation}\label{conj}
c(A,X_1,\dots,X_k)=(AX_1A^{-1},\dots,AX_kA^{-1}),
\end{equation}
and $R(S^3-K)$ is invariant under that action. The quotient can be identified with $X(S^3-K)$.  

Let  
\begin{equation}
C^{i}(W;Ad\rho)=C^i(W)\otimes \g.
\end{equation}  
Each  $C^{i}(W;Ad\rho)$ is an innerproduct space, where the innerproduct is the tensor product of the innnerproduct on 
$C^*(W)$ with the standard dot product on $\mathbb{R}^3=su(2)$. 
There is a chain complex 
\begin{equation}\label{tangential}
\begin{CD} C^0(W;Ad\rho) @>\delta^1>>C^1(W;Ad\rho) @>\delta^2>> C^2(W;Ad\rho) 
\end{CD},
\end{equation}
called the {\em tangential complex}.
Its boundary maps can be understood as the derivatives of the action by conjugation
and the derivative of $\overline{r}$. Its cohomology is closely related to the Zariski tangent space of $X(S^3-K)$,
at a representation $\rho:\pi_1(S^3-K)\rightarrow \G$, as in \cite{W}. 
Call its homology groups the {\em cohomology of $W$ with coefficients in $Ad$ of $\rho$}.

The boundary operators are as follows.  Suppose that $\rho(x_i)=X_i$ and suppose that the relator $r_j$ is a word
in the $x_i$.  First, $\delta^1:C^0(W;Ad\rho)\rightarrow C^1(W;Ad\rho)$ is given by
\begin{equation} 
\delta^1(v^1\otimes \vw)=\sum_i e^i\otimes (X_i-1). \vw. 
\end{equation}
Next, $\delta^2: \oplus e^i\otimes \g \rightarrow \oplus f^j\otimes \g$ is given
by the Fox Jacobian of $(r_1,\ldots,r_{k-1})$ instantiated on $(X_1,\ldots,X_k)$, and
acting by the adjoint action:
\begin{equation}  \delta^2(e^i\otimes \vw)= 
\sum_jf^j\otimes
\frac{\partial r_j}{\partial x_i}.\vw.
\end{equation}
The Fox
calculus and this chain complex are described in \cite{CF,D,J,W}.  Under the standard identification of $T_XSU(2)$ with $\g$ by right translation, $\delta^1$ is the derivative of the conjugation map $c$,  given by (\ref{conj}), at the identity, and $\delta^2$ is the derivative of the relator map 
$\overline{r}$, from equation  (\ref{complex}), at the representation $\rho$.

The tangential complex (\ref{tangential}) has too much homology for its torsion to be interpreted geometrically. To get around that Dubois \cite{DuF} extends the complex so that the torsion defines a one-form on a large portion of the character variety of the knot. 

A representation $\rho:\pi_1(S^3-K)\rightarrow SU(2)$ is {\em regular}, if $Ad\rho$ is irreducible, $H^1(S^3-K;Ad\rho)=\bR$, and the restriction of $\rho$ to the fundamental group of the boundary torus is not central. We will denote regular representations by $R^i(S^3-K)$. The last condition is superfluous for a knot complement since the fundamental group of $S^3-K$ is generated by conjugates of the meridian. Hence the irreducibility of the representation implies that the meridian is not sent to an element of the center of $SU(2)$.  Therefore, $H^0(T^2;Ad\rho)=\bR$ and it is generated by $v^1\otimes \vP$, where
the meridian is sent to $\cos{\phi}+\sin{\phi}\vP$.  Euler characteristic considerations allow us to conclude that if
$\rho$ is regular then $H^2(S^3-K;Ad\rho)=\bR$, and the map 
\begin{equation}\label{firm} 
H^2(S^3-K;Ad\rho)\rightarrow H^2(T^2;Ad\rho),\end{equation}
from the long exact sequence of the pair $(S^3-K,T^2)$
is an isomorphism.

Dubois works with an extended chain complex
\begin{equation}\label{firstsequence} 
\begin{CD} C^0(W;Ad\rho) @>\delta^1>>C^1(W;Ad\rho) @>\delta^2>> C^2(W;Ad\rho) @>\delta^3>> \bR\end{CD}.\end{equation}
To define $\delta^3$ we need to have a firm grasp on the  map (\ref{firm}). It is the adjoint of the inclusion map on homology with coefficients in $Ad\rho$.

Geometrically it is clear that the fundamental class of the boundary torus is homologous to a chain in which
each two cell appears twice, once with positive orientation and once with negative orientation. This chain can be
computed by taking words $L$ and $M$ in the $x_i$ corresponding to the longitude and the meridian and
then writing the word $LML^{-1}M^{-1}$ as a product of conjugates of relators and their inverses.
Each relator $r_j$ will appear twice: once as $s_jr_js_j^{-1}$ and once as $t_jr_j^{-1}t_j^{-1}$. The image
of $[T^2]\otimes \vP$ in $C_2(W;Ad\rho)=C_2(W)\otimes su(2)$ is $ \sum_jf_j\otimes (Ad_{\rho(s_j)}-Ad_{\rho(t_j)}).\vP$. Denote 
\begin{equation}\label{sandt}
\vv_j= (Ad_{\rho(s_j)}-Ad_{\rho(t_j)}).\vP ,
\end{equation}
and
let 
\begin{equation}\label{firstv}
\vv= \sum_jf^j\otimes \vv_j.\end{equation}  
The map $\delta^3:C^2(W;Ad\rho) \rightarrow \bR$ is  the adjoint of the inclusion map of the boundary torus into the knot complement on homology.  It is equal to 
the dot product with $\vv$.

The complex (\ref{firstsequence} ) only has homology in dimension $1$, and that homology at a regular representation $\rho$ is the tangent space at $[\rho]$
to the character variety of the knot complement.

Here is how to pass from a tangent vector to the character variety at $[\rho]$ to a cycle in the chain complex $C^1(W;Ad\rho)$.
For any regular representation $\rho$ there is a neighborhood of $\rho$ in $X(S^3-K)$ which is a smooth $1$-manifold.
Since the quotient map $R(S^3-K)\rightarrow X(S^3-K)$ is a submersion of smooth
manifolds, in a neighborhood of that point we can choose a slice $\rho_t:(-\epsilon,\epsilon) \rightarrow R(S^3-K)$ lifting any smooth
path in $X(S^3-K)$ on a small interval. 
Let $[\rho_t]:[-\epsilon,\epsilon]\rightarrow X(S^3-K)$ be a smooth path.
The cochain  
\begin{equation}
e^i\otimes \left(\frac{d\rho_t(x_i)}{dt}|_{t=0}\right)\rho_0(x_i^{-1})\end{equation}
is a cycle that represents the tangent vector of the path $[\rho_t]$.

Since the chain complex  (\ref{firstsequence}) we are working with has an innerproduct, there is a class of preferred bases used to compute its Reidemeister torsion (\ref{reidtor}).
Let ${\bf{c}}_i$ be an orthonormal basis for the $C^i(W;Ad\rho)$. Let ${\bf{b}}_0={\bf{c}}_0$, and choose 
${\bf{b}}_1$ as in the definition of (\ref{reidtor}).  
Notice that the complement of the space of $2$-cycles is $1$-dimensional. Thus $\bf{b}_2$ consists of a single vector, which we choose to have length $1$.
Recall that the only non-zero homology group is $H^1(W;Ad\rho)$ and that it corresponds to the tangent space of the character variety at $\rho$. Choose a tangent vector 
$h\in T_{[\rho]}X(S^3-K)$.
The absolute value of Dubois' torsion evaluated at $h$ is given by
\begin{equation} 
\tau_{\rho} (h)=
\left|\frac{[\overline{{\bf{b}}}_0,h,{\bf{b_1}}/{\bf{c_1}}]||\vec{v}||}{[\overline{\bf{b}}_1,\bf{b}_2/\bf{c}_2]}\right|. 
\end{equation}
We use the subscript to emphasize the fact that we started with a choice of a regular representation $\rho$.

Dubois works with a signed refinement of torsion, using a scheme that is discussed more generally in Turaev's book \cite{T2}. In fact, Dubois proves that  the torsion defines a $1$-form on the manifold which is the part of the character variety of the knot corresponding to the regular representations. As the regular part of the character variety is oriented there
is a well defined notion of what the sign of the integral of the torsion should be. We can't get the sign from our approach, as it is not inherent to the Gaussian integral that we use to get the global formula.

\subsection{Example continued}

We will compute the torsion at the representations of the fundamental group of the complement of the trefoil
knot along the path $\rho_t$ given by
\begin{equation} \rho_t(x)={\bf{i}}, \ \rho_t(y)=\cos{\frac{\pi}{3}}+\sin{\frac{\pi}{3}}(\cos{t}{\bf{i}}+\sin{t}{\bf{j}}),\end{equation}
where $h$ is the tangent vector of the path. We are using right translation to identify the tangent space at $A$
of $SU(2)$ with $\g$, so that the value of $h$ on the word $w$ in $x$ and $y$ is
\begin{equation} \frac{d\rho}{dt}(w)\rho(w)^{-1}.\end{equation} We are interested specifically in the ordered pair $(h(x),h(y))$ which
is
\begin{equation}(0,\sin{\frac{\pi}{3}}\cos{\frac{\pi}{3}}(-\sin{t}{\bf{i}}+\cos{t}{\bf{j}})+sin^2{\frac{\pi}{3}} {\bf{k}} ).\end{equation}

The boundary operator \begin{equation}\delta^0(v^1\otimes \vv)=e^1\otimes({\bi}-1).\vv +e^2\otimes(\cos{\frac{\pi}{3}}-1+\sin{\frac{\pi}{3}}(\cos{t}{\bf{i}}+\sin{t}{\bf{j}})).\vv\end{equation} can be understood as follows. The coefficient of $e^1$ is $-2$ times the projection onto the perpendicular to $\bf{i}$.
The coefficient of $e^2$ is the projection onto the perpendicular to $\cos{t}{\bi}+\sin{t}{\bj}$ followed by a counterclockwise
rotation by $5\pi/6$ degrees, followed by a homothety of $2\sin{\pi/6}$. Luckily, this will not be that important.
The next boundary operator $\delta^1$ is given by  \begin{equation}\delta^1(e^1\otimes \vv +e^2\otimes \vw)=f^1\otimes ((1+\rho(x)).\vv+(1+\rho(y)+\rho(y)^2).\vw).\end{equation}

We use the following preferred bases for the chain groups.  To start with, $C^0(W;Ad\rho)$ has
basis ${\bf{c}_0}=\{v^1\otimes {\bf{i}}, v^1\otimes {\bf{j}}, v^1\otimes {\bf{k}}\}$.  The basis for $C^1(W;Ad\rho)$ is
 $ {\bf{c}_1}=\{e^1\otimes {\bf{i}}, e^1\otimes {\bf{j}}, e^1\otimes {\bf{k}}, e^2\otimes {\bf{i}}, e^2\otimes {\bf{j}}, e^2\otimes {\bf{k}}\}$. 
 For $C^2(W;Ad\rho)$ we use the basis ${\bf{c}_2}=\{f^1\otimes {\bf{i}},f^1\otimes {\bf{j}}, f^1\otimes{\bf{k}}\}$. Finally for $\bR$ the basis consists of $1$.

As $\rho$ is irreducible, $Z^0(W;Ad\rho)=0$ and we can use the basis above for $C^0(W;Ad\rho)$ to
be $\bf{b}_0$.  A complement of $Z^1(W;Ad\rho)$ is spanned by $e^1\otimes {\bf{i}}, e^2\otimes (cos{t}{\bf{i}} + \sin{t}{\bf{j}})$ so we use this as $\bf{b}_1$. A complement to $Z^2(W;Ad\rho)$ is given by $\bk$ so we use this as $\bf{b}_2$. The last map is onto so the complement of the cocycles in $\bR$ is $0$.  

The change of basis matrix from ${\bf{c}}_0$ to ${\bf{b}}_0$ is the identity so $[{\bf{b}}_0/{\bf{c}}_0]=1$.  The change of basis matrix from $\bf{c}_1$ to $\overline{\bf{b}_0}h\bf{b}_1$ looks like
\begin{equation}\label{bigmatrix} \begin{pmatrix} 0 & 0 & 0 &0 &1 & 0 \\ 0 & -2 &0 & 0 & 0 & 0 \\ 0 & 0 & -2 & 0 & 0 & 0  \\ -\frac{3}{2}+\frac{3}{2}\cos^2(t) & * & * &-\frac{\sqrt{3}}{4}\sin{t} & 0 &\cos{t}  \\\frac{3}{2}\sin{t}\cos{t} & * & *&\frac{\sqrt{3}}{4}\cos{t} & 0 & \sin{t} \\-\frac{\sqrt{3}}{2}\sin{t} & * & * & \frac{3}{4} &0 & 0
  \end{pmatrix}.\end{equation}
  
   The stars are there to emphasize that if you expand the determinant correctly, they do not enter into the computation of the determinant which is $6\sin{t}$.

  Applying $\delta^2$ to $\bf{b}_2$ and adjoining ${\bf{k}}$ we see that the change of basis matrix from ${\bf{c}_2}$ to $\overline{{\bf{b}_2}}{\bf{k}}$ is 
  \begin{equation}\label{littlematrix} \begin{pmatrix} 2 & 3\cos{t} & 0 \\ 0 & 3 \sin{t} & 0 \\ 0 & 0 & 1\end{pmatrix},\end{equation}
  whose determinant is also $6\sin{t}$.
  
  It remains to compute $||\vv||$. The value of the meridian under $\rho$ is,
  \begin{equation} \rho(\mu)=\rho(xy^{-1})=\sin{\frac{\pi}{3}}\cos{t}+\cos{\frac{\pi}{3}}{\bi}-\sin{\frac{\pi}{3}}\sin{t}{\bk}.\end{equation}
  The normalized imaginary part is
  \begin{equation}\vP=\frac{\cos{\frac{\pi}{3}}{\bi}-\sin{\frac{\pi}{3}}\sin{t}{\bk}}{\sqrt{\cos{^2\frac{\pi}{3}}+\sin{^2\frac{\pi}{3}}\sin^2{t}}}.\end{equation}
  Recall that the inclusion map of the chains on the boundary to chains on $W$ takes
  $[T^2]\otimes \vP$ to 
  \begin{equation}
  r\otimes (\rho(x)-\rho(\mu)).\vP=r\otimes \frac{2\sin{\frac{\pi}{3}}\sin{t}{\bk}}{\sqrt{\cos{^2\frac{\pi}{3}}+\sin{^2\frac{\pi}{3}}\sin^2{t}}}.
  \end{equation}
  
  Putting it all together 
  \begin{equation}\label{tortref}
  \tau (h)=\frac{2\sin{\frac{\pi}{3}}\sin{t}}{\sqrt{\cos{^2\frac{\pi}{3}}+\sin{^2\frac{\pi}{3}}\sin^2{t}}}.
  \end{equation}
  Note that this agrees with the Proposition 5.4 in \cite{DuF}, after you account for differences
  in notation, and parametrization.

Let's interpret this computation geometrically. The determinant of the matrix  (\ref{bigmatrix}) computes the volume form on $su(2)\oplus su(2)=C^1(W,Ad\rho)$ which is the
Riemannian volume form on $T_{(X,Y)}SU(2)\times SU(2)$ at the ordered pair corresponding to the representation $\rho_t$. The last two columns are an orthogonal complement to the image
of $T_{\rho}R(S^3-K)$ at that point. Up to sign the determinant of this matrix
is computing the Riemannian volume form on the representation variety of the knot at $\rho$ inherited from its embedding in $SU(2)\times SU(2)$.  This means that we can interpret
the determinant of matrix (\ref{bigmatrix})  as  
\begin{equation}
\nu^R\begin{pmatrix} 0 & 0 & 0 &0  \\ 0 & -2 &0 & 0  \\ 0 & 0 & -2 & 0  \\ -\frac{3}{2}+\frac{3}{2}\cos^2(t) & * & * &-\frac{\sqrt{3}}{4}\sin{t} \\ \frac{3}{2}\sin{t}\cos{t} & * & *&\frac{\sqrt{3}}{4}\cos{t} \\-\frac{\sqrt{3}}{2}\sin{t} & * & * & \frac{3}{4} 
  \end{pmatrix}\end{equation}
where $\nu^R$ is the Riemannian volume form on the representation variety of the knot evaluated at the column vectors of the matrix.
Further contemplation of the formula reveals that the argument of $\nu^R$ is the
derivative of $\beta(A,t)$ at $A=Id$, where $\beta$ is the parametrization given by equation (\ref{beta}).  The value of the determinant of (\ref{bigmatrix}) is then equal to
\begin{equation}
(\beta_{(Id,t)})^*
(\nu^R)
({\bf{i}},{\bf{j}},{\bf{k}},\frac{d}{dt}).
\end{equation}
The formula above means that we have used $\beta$ to pull the volume form $\nu^R$ back
to $T_{(Id,t)}SU(2)\times (0,\pi)$ and then evaluated that pulled back form on an orthonormal basis.

The first three columns of the matrix (\ref{bigmatrix}) make up a basis for the tangent space of the fiber of the projection map from the representation variety of the knot complement to the character variety at $[\rho_t]$. If we had chosen an orthonormal basis for this space, the determinant would be computing the push-forward volume on the character variety. However, the basis we chose is not orthonormal.  Letting $c:SU(2)\rightarrow SU(2)^2$ be the result of fixing $t$ in $\beta$, that is  $c(A)=\beta(A ,t)$,
we could interpret this formula for fixed $t$ as
\begin{equation} ||dc_{Id}||(\beta_{(Id,t)})^*\nu^X(\frac{d}{dt}),\end{equation}
where $\nu^X$ is the Riemannian volume form from the push-froward metric on the character variety of $\pi_1(S^3-K)$
at the point $[\rho_t]$.
Notice that $dc_{Id}$ is the first boundary operator in the complex used to define the Dubois torsion. 

The determinant  of the matrix  (\ref{littlematrix}) is  $\pm ||d\overline{r}||$. Finally,
the last factor in the formula for torsion contributes $||\vec{v}||$.

Putting it all together this means
that the Reidemeister torsion can be thought of as a one-form $d\tau$ on the character variety of the knot complement, given by
\begin{equation} d\tau=\frac{||dc||\ ||\vec{v}||}{||d\overline{r}||}\nu^X,\end{equation}
where $dc$ is the derivative of the map given by conjugation, $d\overline{r}$ is the derivative of the relator map, and $\nu^X$ is the push-forward of the Riemannian volume from the representation variety to the character variety.  

However, a more useful formulation for us takes place in the representation variety, where
we see

\begin{equation} d\tau(h)=\frac{\nu^R(D\beta_{(Id,t)})||\vec{v}||}{||d\overline{r}||},\end{equation}
where $\beta:SU(2)\times (-\epsilon,\epsilon)\rightarrow R^i(S^3-K)$ is of the form
$\beta(A,t)=A\rho_tA^{-1}$ and $\rho_t:(0,\pi)\rightarrow R^i(S^3-K)$ is a slice
of the projection mapping $R^i(S^3-K)\rightarrow X(S^3-K)$ over the path $[\rho_t]:(0,\pi)\rightarrow X(S^3-K)$ whose derivative is the tangent vector $h$. Finally, $\nu^R$ is the
Riemannian volume form on $R^i(S^3-K)$ coming from its inclusion in $SU(2)^k$.

  \subsection{Construction of the seminorm}

 A {\em seminorm} on a vector space $V$ is a map 
$||\  ||_s:V \rightarrow \bR_{\geq 0}$ such that for all $\vv,\vw \in V$,
$||\vv+\vw||_s\leq ||\vv||_s+||\vw||_s$, and for any $\lambda\in \bR$, 
$||\lambda \vv||_s=|\lambda|||\vv||_s$.  If $V$ is in addition
a commutative algebra, we define the {\em radical} of $|| \ ||_s$ to be the set of all 
$\vv \in V$ such that for all $\vw \in V$,
$||\vv \vw||_s=0$.  By design the radical of a norm is an ideal. Any norm is a seminorm and its radical is the zero ideal.

At a regular representation, the sequence
\begin{equation}\label{secondseq} 
\begin{CD}0@>>>T_{\rho}R(S^3-K)@>i>>\oplus_k su(2)@>\delta^2>> \oplus_{k-1} su(2) @>\delta^3>> \bR@>>>0\end{CD}
\end{equation}
is exact. 
This is almost the same sequence as (\ref{firstsequence}). The sequences differ in the first term, and the  map $i$ is an inclusion. 
We use different notation in the second and third terms. For instance,
$\sum_ie^i\otimes w_i\in C^1(W, Ad\rho)$ corresponds to the tuple $(w_1,\ldots w_k)\in \oplus_k su(2)$.
Using direct sums of the standard orientations on $\bR$ and on $\g=\bR^3$, we have orientations for $\oplus_k\g$ and $\oplus_{k-1}\g$. 
This allows
us to force an orientation on $T_{\rho}R(S^3-K)$. As $T_{\rho}R(S^3-K)$ inherits an innerproduct as
a subspace of $\oplus_k\g=\bR^{3k}$ there is a Riemannian volume $\nu^R$ on $T_{\rho}R(S^3-K)$.

The torsion of (\ref{secondseq}) defines a $4$-form on $T_{\rho}R(S^3-K)$. Using the action of $SU(2)$ on $R(S^3-K)$ by conjugation we can push the $4$-form down to a $1$-form on $T_{[\rho]}X(S^3-K)$.  The push-forward is Dubois' torsion.

In our example, we could parametrize the regular representations of the fundamental group of the complement of the trefoil by taking a path of representations and conjugating them.  You can do this at any regular representation of a knot group. Here is how: Parametrize $X(S^3-K)$ near $[\rho]$  by 
$[\rho_t]:(-\epsilon,\epsilon)\rightarrow X(S^3-K)$.
Using a slice  of
the quotient map $R(S^3-K)\rightarrow X(S^3-K)$ lift $[\rho_t]$ to
$\rho_t:(-\epsilon,\epsilon)\rightarrow R(S^3-K)$.
Combining with conjugation we 
get a local parametrization $\beta:SU(2)\times (-\epsilon,\epsilon)\rightarrow R(S^3-K),$
\begin{equation}\label{paramet} 
\beta(A,t)=A\rho_tA^{-1}. 
\end{equation}

Computing $\nu^R$ in these coordinates we get that
it is $[\overline{\bf{b}}_0,\rho'_t\rho_t^{-1},\bf{b}_1/c_1]$ up to sign, as long as we chose
our ${\bf{b}}_0$ and ${\bf{b}}_1$ to be orthonormal.  Furthermore, if we choose $\bf{b}_2$
to be orthonormal as well then $[\overline{\bf{b}}_1,\overline{\bf{b}}_2/{\bf{c}}_2]$ is 
$||d{\overline{r}}||$ up to sign, where ${\overline{r}}$ was defined in (\ref{complex}).
Thus
if we orient $R^i(S^3-K)$  as above then the absolute value of Dubois' torsion at
a representation $\rho$
is  
\begin{equation}
\tau_{\rho}(h)= 
\frac{\nu^R(D\beta_{(Id,t)})||\vec{v}||}{||d{\overline{r}}||}.
\end{equation}
The manifold $R^i(S^3-K)$ is oriented, and the local parametrizations we have chosen
preserve orientation. Hence we can integrate  torsion in local coordinates. By standard arguments these integrals can be woven together to yield a global
value so long as the value of the integral is finite.
Hence, for any
function  $f:R^i(S^3-K)\rightarrow \mathbb{R}$ its integral against Dubois' torsion is given by
\begin{equation}\label{inttor}
\int_{R^i(S^3-K)}fd\tau=
\int_{R^i(S^3-K)} f \frac{||\vec{v}||}{||d{\overline{r}}||}\nu^R,
\end{equation}
provided that it converges. 

\begin{defi}  Let $f\in C[X(T^2)]$, using restriction we can view $f$ as a function on $R(S^3-K)$. Define
\begin{equation}\label{seminorm}
||f||_s=\left|\int_{R^i(S^3-K)} f \frac{||\vec{v}||}{||d\overline{r}||}\nu^R\right|  \end{equation}
if the integral is defined.\end{defi}

\begin{theorem}
If the regular representations $R^i(S^3-K)$ form an open dense subset of $R(S^3-K)$ and the total Dubois' torsion
of $R^i(S^3-K)$ is finite,
then the seminorm $||\ ||_s$ is a knot invariant whose radical is the $A$-ideal. 
\end{theorem}
\proof
This is follows from the fact that when it is defined,
Reidemeister torsion is never $0$, so the support of the measure defined by the Dubois' torsion is the closure of the  regular irreducible representations.\qed

Define the {\em total Dubois torsion of the knot} $K$ to be
\begin{equation}\label{totaltor}
\tau(K)= \int_{X^i(S^3-K)}d\tau.
\end{equation}
 In the case of torus knots, the total torsion is exactly computable.  For instance the total
Dubois torsion  of the trefoil knot
is $\frac{4\pi}{3}$.  It is a little more work, but completely elementary to compute the integral of any peripheral character against Dubois torsion.  Calculate the integrals of the powers of the real part of the image of the meridian, and then prove these are a basis for the values of the integral against Dubois torsion of the peripheral skeins.

\section{A global computation of the seminorm}
The goal of this section is to derive a global formula for integrating against Reidemeister torsion. We introduce a function on $SU(2)^k$ that involves a parametrix of the  heat kernel trace, and show that the limit of integrals against that function yields the seminorm defined above. The proof involves an auxiliary family of integrals to which Laplace's method can be applied. Finally, we replace the parametrix of the heat kernel with  the heat kernel trace, to obtain a global formula.

\subsection{Two pointwise close Dirac delta functions}

We are working with  Cartesian products of $SU(2)$ and intervals in the real line.
We treat $SU(2)$ as the unit sphere in $\bR^4$. The volume form $\nu^{SU(2)}$ from the Riemannian metric is $2\pi^2$ times
Haar measure.  A good parametrization of $SU(2)$ is given by three angles $(\phi,\psi,\theta)$, where $\phi$ and $\psi$ vary from $0$ to $\pi$, and $\theta$ varies from $0$ to $2\pi$.  The parametrization is
\begin{equation}
X(\phi,\psi,\theta)=(\cos{\phi},\sin{\phi}\sin{\psi}\cos{\theta},\sin{\phi}\sin{\psi}\sin{\theta},\sin{\phi}\cos{\psi}).
\end{equation}
The angle $\phi$ is the angle that $X(\phi,\psi,\theta)$ makes with $(1,0,0,0)={\bf 1}$.
The projection of $X(\phi,\psi,\theta)$ into the perpendicular to $(1,0,0,0)$ is a vector in $\mathbb{R}^3$ with length $\sin{\phi}$ that makes an angle $\psi$ with the $z$-axis, and whose projection into the $xy$-plane makes an angle $\theta$ with positive $x$-axis. The volume form  from the Riemannian metric inherited from $\bR^4$ is 
\begin{equation}
\nu^{SU(2)}=
\sin^2{\phi}\sin{\psi} \ d\phi\wedge \ d\psi \wedge\ d\theta .
\end{equation}

There are two approximate Dirac delta functions of ${\bf 1}$ with respect to the measure $\nu^{SU(2)}$ we would like to use.
The first is
\begin{equation}
\beta_{\lambda}(\phi,\psi,\theta)= \left(\frac{\lambda}{4\pi}\right)^{3/2}e^{-\lambda\phi^2/4}.
\end{equation}

There is a very nice coordinate system for working with this Gaussian.  Recall the exponential map,
$exp:\mathbb{R}^3\rightarrow S^3$, given by 
\begin{equation}
exp(\vw)=\cos{||\vw||}+\sin{||\vw||}\frac{\vw}{||\vw||}.  
\end{equation}
When restricted to the open ball of radius $\pi$ centered at the origin,  $B_{\pi}(\vec{0})$, the exponential map is a diffeomorphism onto $S^3-\{-{\mathbf{1}}\}$, and has an inverse
\begin{equation}\log: S^3-\{{-\mathbf{1}}\}\rightarrow B_{\pi}(\vec{0}),\end{equation}
which is a normal coordinate chart.  In these coordinates,
\begin{equation}\beta_{\lambda}(\vw)=\left(\frac{\lambda}{4\pi}\right)^{3/2}e^{-\lambda||\vw||^2/4}.\end{equation}

Suppose that $r:M\rightarrow SU(2)$ is a smooth function from a smooth manifold $M$ so that $r(p)={\mathbf{1}}$, and let
$U$ be a coordinate chart with coordinates $x^i$ at $p$, so that $w=log(r)$ is defined
in $U$.  Since $exp(w)=r$ we have that for any $\frac{\partial}{\partial x^i}|_p$,
\begin{equation}  \frac{\partial}{\partial x^i}r|_p=\frac{\partial}{\partial x^i}exp(w)|_p=exp(w(p)) \frac{\partial}{\partial x^i}w.\end{equation}
 Since $exp(w(p))=r(p)={\mathbf{1}}$, we have that
 \begin{equation}\label{logs}
 Dr|_p=Dw|_p.\end{equation}

The second approximate Dirac delta function is the {\em heat kernel trace} $\kappa_{\lambda}$. The reason for working with $SU(2)$ as the $3$-sphere in ${\mathbb{R}}^4$ is to have the eigenvalues of the Laplace operator be integers.
 In the 
$(\phi,\psi,\theta)$ coordinates,
\begin{equation} \kappa_{\lambda}(\phi,\psi,\theta)=\frac{1}{2\pi^2}\sum_{c=0}^{\infty}(-1)^{c}(c+1)e^{-c(c+2)/\lambda}s_c(-2\cos{\phi}),\end{equation}
where $s_c$ is the $c$th Chebyshev polynomial defined by the recursion,
$s_0=1$, $s_1=x$, and $s_n=xs_{n-1}-s_{n-2}$.  In fact, $s_c(-2\cos{\phi})$ is equal to 
$(-1)^ctr_c$ where $tr_c$ is the trace in the $(c+1)$-dimensional irreducible representation of $SU(2)$. These Chebyshev polynomials are the same as the ones used to define the colored Jones polynomial.

The relationship between these two approximate Dirac delta functions is that the first is the parametrix for the heat kernel trace and the second is the heat kernel trace.  By design they are pointwise close. 
\begin{lemma}\label{Migdal}
  For every $\epsilon>0$ there exists $C>0$ and $U$ a neighborhood of ${\bf{1}}$ in $SU(2)$ so that
for $\lambda \geq C$, and $A\in U$,
\begin{equation} |\beta_{\lambda}(A)-\kappa_{\lambda}( A)|<\epsilon.\end{equation}\end{lemma}

\proof This follows directly from lemma 3.18 of \cite{R}, by substituting $\lambda$ for $1/t$. \qed

The fact that these two Dirac delta functions are asymptotically equivalent is sometimes called {\em Migdal's equivalence} \cite{Mi}.  We were first exposed to this type of local-global argument by reading \cite{BuRo}. Recently Bozom and Smerlak used the same approach to understand the Reidemeister torsion of two complexes \cite{BoSm}.

\subsection{Laplace's method}

Let $M$ be a Riemannian manifold with Riemannian volume form $\nu^M$. Let
$\Phi:M\rightarrow \mathbb{R}$ be a non-negative function that takes on the value $0$ along the codimension $k$ submanifold $N$ with induced Riemannian volume form $\nu^N$.
Let $H_N(\Phi)$ be the restriction of the Hessian of $\Phi$ to the normal space
to $N$.  Assume that $H_N(\Phi)$ is nondegenerate at each point in $N$ and there
is an $\epsilon$-neighborhood of $N$ outside of which the function $\Phi$ is bounded away from $0$. It is an easy computation in local coordinates to see that if the integral on the
left exists for the smooth function $f:M\rightarrow \mathbb{R}$, then
\begin{equation}
\label{Laplace} \lim_{\lambda\rightarrow \infty} \left(\frac{\lambda}{\pi}\right)^{k/2}\int_Mfe^{-\lambda \Phi}\nu^M=\int_N\frac{f\nu^N}{\sqrt{(det(H_N))}}.\end{equation}
The equation is sometimes referred to as {\em Laplace's method} \cite{Wo,BH}.

We will need to use Laplace's method  in a slightly more general context.  
Given  a smooth non-negative real-valued function $\Phi$   defined on a subset of $M$,  let $N=\{p\in M \ | \ \Phi(p)=0\}$.
Let $S$ be the subset of $M$ consisting of three types of points:
\begin{itemize}
\item Points where $\Phi$ is not defined;
\item Points where $N$ is not a manifold;
\item Points where $N$ is a manifold but the Hessian of $\Phi$ restricted to the normal direction to $N$ is degenerate.
\end{itemize}
Let $U_{\epsilon}$ be an epsilon neighborhood of $S$.  If the limit 
\begin{equation} 
C(\epsilon)=\lim_{\lambda\rightarrow \infty}\left(\frac{\lambda}{\pi}\right)^{k/2}\int_{U_{\epsilon}}e^{-\lambda \Phi}\nu^M\end{equation} 
exists, and $\lim_{\epsilon\rightarrow 0}C(\epsilon)=0$ then we say the singularities of $\Phi$ are {\em tame}.

\begin{prop}  Let $M$ be compact, and suppose that the singularities $S$ of $\Phi$ are tame,
and $f:M\rightarrow \mathbb{R}$ is smooth, then
\begin{equation}\label{laplace} \lim_{\lambda\rightarrow \infty} \left(\frac{\lambda}{\pi}\right)^{k/2}\int_Mfe^{-\lambda \Phi}\nu^M=\int_{N-S}\frac{f\nu^N}{\sqrt{(det(H_N))}}.\end{equation}
\end{prop}   \qed

\subsection{The Integral}
We will define a function whose maximum occurs along representations of the knot complement. We use Laplace's method to compute an appropriately normalized limit of an integral of this function. The answer turns out to be the seminorm defined by Doubois torsion. This function depends on additional parameter, which can be integrated out to recognize the parametrix of the heat kernel. 

Recall that given a knot $K\subset S^3$, we denoted the regular representations of $\pi_1(S^3-K)$ into $SU(2)$ by $R^i(S^3-K)$ and the abelian
 representations of $\pi_1(S^3-K)$ into $SU(2)$ by $Ab(S^3-K)$.
We work with a $2$-complex embedded in the complement of the knot $K$ as a deformation retract,
giving rise to $k-1$ relator maps corresponding to the two-cells, $(r_1,\ldots,r_{k-1})$.  The set $S\subset SU(2)^k$,
which is the locus along which any of the relators takes on the value $-1$,  has measure zero as long as there
is  some regular representation at which the variety $R^i(S^3-K)$ has dimension $4$. For a knot in $S^3$ this is equivalent to the statement that the set of regular representations $R^i(S^3-K)$ is nonempty.

Notice that $\vec{v}$ defined in (\ref{firstv}) can be described by word maps that are well defined for all points in $SU(2)^k$. 
Let
 \begin{equation}
 \overline{v}=\frac{1}{|| (\vv_1,\ldots,\vv_{k-1})||^2}
 \sum_j f^j\otimes \vv_j.
\end{equation}
The normalized vector $\overline{v}$  is not defined on a set consisting of points where $\vec{v}=\vec{0}$. Call this set of points $T$. Thus  the vector  $\overline{v}$   can be extended to a function  $\overline{v}:SU(2)^k - T\rightarrow su(2)^{k-1}$. 
Once again if the set of regular representations of the fundamental group of the knot is nonempty then $T$ has measure $0$.
 
Notice that $R^i(S^3-K)\cap T=\emptyset$  and $Ab(S^3-K)\subset T$.  This is because 
 
 \begin{equation}  
 \vv=((Ad_{\rho(s_1)}-Ad_{\rho(t_1)}).\vec{P},\ldots, (Ad_{\rho(s_{k-1})}-Ad_{\rho(t_{k-1})}).\vec{P}),
 \end{equation}  
 where $s_i,t_i$ where defined by (\ref{sandt}).
 At abelian representations the vector $\vec{P}=\frac{Im\rho(\mu)}{||Im\rho(\mu)||}$, where $\mu$ is the meridian, is fixed
 by $Ad_{\rho(\gamma)}$ for all $\gamma\in \pi_1(S^3-K)$. Thus the vector $\vv=\vec{0}$,
 and its norm is zero.
 
  We will use the tuple notation as opposed to the tensor used in (\ref{firstv}).
 \begin{equation}\label{normv}
 ||\vv||=||(\vv_1,\ldots\vv_{k-1})||,
 \end{equation}
 and
 \begin{equation}\label{newv}
 \overline{v}=\frac{1}{||\vv||^2}(\vv_1,\ldots\vv_{k-1})= (\overline{v_1},\dots , \overline{v}_{k-1}).
 \end{equation}

  Fixing a relator $r_i$, let
$\vw_i=\log{r_i}$.   We have that all $\vw_i$ are well defined as functions from $SU(2)^k-S$ to $\mathbb{R}^3$.  As each
$\overline{v}_i$ is well defined as a function from $SU(2)^k-T$ to $\mathbb{R}^3$ we can put it all together to define
\begin{equation}\label{Phi}
\Phi=\frac{\sum_i||\vw_i||^2+t^2||\overline{v}_i||^2}{4}:\left(SU(2)^k-(S\cup T)\right)\times \mathbb{R}\rightarrow \mathbb{R}.
\end{equation}
This function attains its global minimum along $R^i(S^3-K)\times \{0\}$, and that minimum is equal to $0$.

\begin{theorem}[Local Computation]\label{Local}   Suppose that the map $\Phi$ defined by 
(\ref{Phi}) has tame singularities. Let $f:SU(2)^k \rightarrow \mathbb{R}$ be any smooth invariant function. Then
\begin{equation}\label{local}
\int_{R^i(S^3-K)}fd\tau=
 \lim_{\lambda\rightarrow \infty}\sqrt{\frac{4\pi}{\lambda}} \int_{SU(2)^k}2^{(3k-3)/2} f \prod_j\beta_{\lambda}(r_j)||\vec{v}|| \nu^k ,\end{equation}
where 
$\nu^k$ is the volume form from the inclusion of $SU(2)^k$ into $\mathbb{R}^{4k}$, and $d\tau$ is the measure from Dubois torsion.
\end{theorem}

\proof Let $f:SU(2)^k\rightarrow \mathbb{R}$ be a smooth function. Consider the integral,

\begin{equation}
\left(\frac{\lambda}{4\pi}\right)^{(3k-3)/2}\int_{SU(2)^k\times \mathbb{R}}f e^{-\lambda \Phi}\nu^k d\mathcal{L},
\end{equation}
where  $d\mathcal{L}$ is the Lebesgue measure on the real line.

If the singularities of $\Phi$ are tame, then we can evaluate the limit via Laplace's method.
Hence, we need to compute its Hessian along $R^i(K)\times\{0\}$.  Notice that $\Phi$ involves a sum of terms of the form $w\cdot w/4$ where $w$ is a function from a Riemannian manifold into
an innerproduct space, and the critical values occur on the locus $w=0$.  Computing the derivative of a term of
this form using the product rule we get $D(w\cdot w)/4=(1/2)Dw\cdot w$. Differentiating again we get,
\begin{equation} (1/2)(D^2w\cdot w+Dw\cdot Dw).\end{equation}  This handles the Hessian of the terms involving $||\vw_i||^2$.  The terms
of the form $t^2||\overline{v}_i||^2$ are easier as we are only interested  along the locus $t=0$. Hence the only contributions come from differentiating with respect to $t$ twice.

Since we  only care about the case  when $w=0$ we get 
\begin{equation}H(w)=(1/2)Dw\cdot Dw.\end{equation}  
Since the sum of the Hessians is the Hessian of the sum we have that along $R^i(S^3-K)\times \{0\}$,
\begin{equation}
H\left(\frac{\sum_i||\vw_i||^2+t^2||\overline{v}_i||^2}{4}\right)=(1/2)\sum_i \begin{pmatrix}  D\vw_i\cdot D\vw_i & 0 \\ 0 &\overline{v}_i\cdot \overline{v}_i\end{pmatrix},\end{equation}
 where we split
the tangent space of $\left(SU(2)^k-(S\cup T)\right)\times \mathbb{R}$ at $(\rho,0)$ as
\begin{equation} T_{\rho}SU(2)^k \oplus T_0\mathbb{R}.\end{equation}   

From the properties of the logarithm, and the fact that we are working along the locus
where the relators take on the value $1$, we have that $D\vw_i=Dr_i$ (compare with (\ref{logs})).  The Hessian splits
as a direct sum over orthogonal matrices so that 
\begin{equation}\label{Hessian} 
\sqrt{det(H_N(\Phi))}=\frac{1}{2^{(3k-3)/2}}\frac{||d\overline{r}||}{||\vv||}.
\end{equation}

Assuming that the singularities of $\Phi$ are tame
with respect to the Gaussian, we see that $S\cup T$ has measure zero. The functions $f$ are bounded, so we can consider the following limit.

\begin{equation}\label{limint}
 \lim_{\lambda\rightarrow \infty}\left(\frac{\lambda}{4\pi}\right)^{(3k-3)/2}\int_{SU(2)^k\times \mathbb{R}} f \ exp\left(-\lambda\frac{\sum_i||\vw_i||^2+t^2||\overline{v}_i||^2}{4}\right)\nu^kd\mathcal{L}.
 \end{equation}
By Laplace's method this is equal to
\begin{equation}
\left(\frac{1}{4}\right)^{(3k-3)/2}\int_{R^i(S^3-K)}\frac{ f\nu^R}{\sqrt{det(H_N(\Phi))}},
\end{equation}
where $\nu^R$ is the Riemannian metric on $R^i(S^3-K)$ inherited from $SU(2)^k$ and $|det(H_N(\Phi)|$ is the determinant of the Hessian of $\Phi$ restricted to the normal of the tangent space to $R^i(S^3-K)$ in $SU(2)^k$.
Substituting the value for the Hessian computed in (\ref{Hessian}) this is further equal to
\begin{equation}
\left(\frac{1}{4}\right)^{(3k-3)/2}2^{(3k-3)/2} 
\int_{R^i(S^3-K)}
\frac{ f  {||\vv||} \nu^R} { ||d\overline{r}|| }
\end{equation} 
 By (\ref{inttor}) this is a multiple of the integral against the Dubois torsion,  
 \begin{equation}\frac{1}{2^{(3k-2)/2}}\int_{R^i(S^3-K)}fd\tau.\end{equation}

Going back to (\ref{limint}), the dependence on $t$ can be factored out so that it is of the form
\begin{equation}
\left(\frac{\lambda}{4\pi}\right)^{(3k-3)/2}\int_{SU(2)^k\times \mathbb{R}} f 
\ exp\left(-\lambda\frac{\sum_i||\vw_i||^2+t^2||\overline{v}_i||^2}{4}\right) \nu^kd\mathcal{L}=
\end{equation}
\begin{equation}\left(\frac{\lambda}{4\pi}\right)^{(3k-3)/2}\int_{SU(2)^k} f \ exp\left(-\lambda\frac{\sum_i||\vw_i||^2}{4}\right)\left(\int_{\mathbb{R}}exp(-\lambda \frac{t^2||\vec{v}||}{4})d\mathcal{L}\right)\nu^k.
\end{equation}

Next we integrate out the $t$ to get,
\begin{equation}\sqrt{\frac{4\pi}{\lambda}}\left(\frac{\lambda}{4\pi}\right)^{(3k-3)/2}\int_{SU(2)^k} f \ exp\left(-\lambda\frac{\sum_i||\vw_i||^2}{4}\right)||\vec{v}||\nu^k.\end{equation}

Putting this into the equation derived from Laplace's method we have,
\begin{equation}\lim_{\lambda\rightarrow \infty}\sqrt{\frac{4\pi}{\lambda}}\left(\frac{\lambda}{4\pi}\right)^{(3k-3)/2}\int_{SU(2)^k} f 
\ exp\left(\lambda\frac{\sum_i||\vw_i||^2}{4}\right)||\vec{v}||\nu^k=
\frac{1}{2^{(3k-2)/2}}
\int_{R^i(S^3-K)}fd\tau.
\end{equation}

Finally, we recognize 
\begin{equation}\left(\frac{\lambda}{4\pi}\right)^{(3k-3)/2} 
exp\left(-\lambda\frac{\sum_i||\vw_i||^2}{4}\right)=\prod_j\beta_{\lambda}(r_j),\end{equation} 
which ends the proof of the  theorem. \qed

Since the Gaussian is pointwise close to the heat kernel trace
we can replace $\beta_{\lambda}$  by $\kappa_{\lambda}$ in the right hand side of the equation (\ref{local})  to get
\begin{equation} \lim_{\lambda\rightarrow \infty}  \sqrt{\frac{4\pi}{\lambda}}\int_{SU(2)^k}
2^{\frac{3k-3}{2}}
f||\vv||\prod_{i=1}^{k-1} \kappa_{\lambda}(r_i)\nu^k,\end{equation}
which yields a global formula for Dubois torsion.

\begin{theorem}[Global Formula]\label{global}
\begin{equation}\int_{X^i(S^3-K)}fd\tau=\lim_{\lambda\rightarrow \infty} 2^{\frac{3k-3}{2}} \sqrt{\frac{4\pi}{\lambda}}\int_{SU(2)^k}\
f||\vv||\prod_{i=1}^{k-1} \kappa_{\lambda}(r_i)\nu^k\end{equation}
\end{theorem}

\proof This follows directly from Theorem \ref{Local} and Lemma \ref{Migdal}. 

\section{ Interpreting the global formula}
The purpose of this section is to interpret integration against Dubois torsion in terms of quantum invariants of the knot complement.

By theorem \ref{global},
\begin{equation}\label{globalformula}
\int_{X^i(S^3-K)}fd\tau =\end{equation} \begin{equation}lim_{\lambda\rightarrow \infty} 2^{\frac{3k-3}{2}} \sqrt{\frac{4\pi}{\lambda}}\int_{SU(2)^k}
f||\vv||\prod_{i=1}^{k-1}
\frac{1}{2\pi^2}\sum_{c=0}^{\infty}(-1)^{c}(c+1)e^{-c(c+2)/\lambda}s_c(-tr (r_i))\nu^k.\end{equation}
In order to relate this formula to quantum invariants we need to recall the definition of the Yang-Mills measure \cite{BFK1}. 

\subsection{The Yang-Mills measure in a handlebody}

In this section we recall definitions of the Kauffman bracket skein module of a handlebody and of the Yang-Mills measure functional on that module.  A good reference for standard definitions  is the book \cite{Li}. More details and the proofs relating to the Yang-Mills measure can be found in \cite{BFK1}.

Recall that the Kauffman bracket skein module of a manifold is defined as a quotient of a vector space over ${\mathbb{C}}$, with basis given by the set of equivalence classes of framed links in the manifold (including an empty link), by the relations that define the Kauffman bracket:
\begin{equation}
\bigcirc\cup L= -(t^2+t^{-2})L
\end{equation}
and
\begin{equation}
\raisebox{-4pt}{\includegraphics{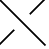}}=t \raisebox{-3pt}{\includegraphics{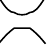}}+t^{-1}\raisebox{-3pt}{\includegraphics{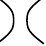}}.
\end{equation}
The elements of  the Kauffman bracket skein module are called skeins. 

Let $H$ be a handlebody.
There is a convention of modeling skeins on admissibly colored framed trivalent graphs in $H$. An admissible coloring is an assignment of a nonnegative integer to each edge, so that the colors at trivalent vertices satisfy all possible triangle inequalities. The skein corresponding to such a graph is obtained by inserting an appropriate Jones-Wenzl idempotent at each edge, inserting a Kauffman triad at each vertex, and joining those with parallel strands. For the definitions of the Jones-Wenzl idempotents and the Kauffman triads see \cite{Li}.

The Yang-Mills measure is a local, diffeomorphism invariant trace defined on $K_t(H)$, 
\begin{equation}
{\mathcal{YM}} : K_t(H) \rightarrow {\mathbb{C}}.
\end{equation}
Here is 
how to compute  ${\mathcal{YM}}$. Given a handlebody $H$ of genus $g$ its double is diffeomorphic to the connected sum of $g$ copies of $S^1\times S^2$. The Yang-Mills measure ${\mathcal{YM}}$ is the linear functional given by taking the inclusion of $H$ into 
$\sharp_gS^1\times S^2$ followed by the canonical isomorphism of $
K_t(\sharp_gS^1\times S^2)$ with $\mathbb{C}$. This isomorphism can be described in the following way: Choose a system of spheres that cut down $\sharp_g S^1\times S^2$ to a punctured ball. Represent a skein as a linear combination of colored trivalent graphs intersecting the spheres transversely in the interior of edges, with each graph intersecting any sphere at most once. If a graph intersects a non-separating sphere in a single point on one edge then this graph represents zero in the skein module.
Thus  we can assume that the graphs miss the spheres. The Yang-Mills measure is the value of the Kauffman bracket of the resulting skein in the punctured ball. 

Alternatively, a handlebody $H$  can be thought of as $H=F\times I$, where $F$ is a compact oriented surface with boundary. Choosing a trivalent spine for $F$ yields a basis for $K_t(F)=K_t(H)$ given by the skeins corresponding to all possible admissible colorings of that spine. The Yang-Mills measure of any skein is the coefficient of the skein coming from labeling all the edges of the spine with $0$. Note that this does not apply to the skein algebra of the disk (which is spanned by the empty skein) nor to the skein module of the annulus (i.e. of the solid torus), which is spanned by all nonnegative colorings of the core of the annulus.

When $t=e^{\frac{\pi i}{2r}}$, the Yang-Mills measure is defined the same way as above on the {\em reduced} Kauffman bracket skein module, $K_{r,f}(H)$, which is obtained from $K_t(H)$ by taking its quotient by the submodule spanned by all the skeins corresponding to the framed trivalent graphs where some edge is colored with $r-1$. The canonical basis for $K_{r,f}(H)$ is finite  since colors cannot be larger than $r-2$.
In this case the admissibility condition carries an additional requirement that the sum of the  three colors at any vertex  is less than or equal to $2r-4$. 

The projector $\omega$ is an element of 
the skein algebra of the solid torus defined below:
\begin{equation}\label{omega}§
\omega =\sum_{i=0}^{r-2} (-1)^i[i+1]s_i,
\end{equation}
where $[n]$ denotes the quantum integer,
\begin{equation}
[n]=\frac{t^{2n}-t^{-2n}}{t^2-t^{-2}},
\end{equation}
and $s_i$ denotes the skein in the annulus which is the result of coloring the core with the $i$-th Jones-Wenzl idempotent. Note that $[n]$ is defined for all $t\neq 0$. At $t=\pm1$ the formula we gave has a removable singularity and $\lim_{t\rightarrow \pm1}[n]=n$. As the level $r$ approaches infinity, $t$ approaches $-1$.

Evaluating the Yang-Mills measure of a skein $s$ in the handlebody $H$ can be understood in terms of
the Kauffman bracket of a surgery diagram in the three-sphere.  Embed $H$ in
$S^3$ so that its complement is also a handlebody. Put an unknotted  zero-framed circle decorated with the projector $\omega$ around each handle of $H$, and then take the value of the Kauffman bracket in the $3$-ball of the skein which is a union of those decorated circles with the skein $s$. Finally
divide the result by $X^{2g}$, where 
\begin{equation}\label{defofX}
X=\sqrt{\sum_{c=0}^{r-2}[c+1]^2}.
\end{equation}
The fact that this agrees with the definition of the Yang-Mills measure follows from the Turaev-Wenzl identity, fusion, and the shadow world formula for the Yang-Mills measure \cite{FKB}.

When the parameter $t$ is equal to $-1$, the algebra $K_{-1}(H)$ is isomorphic to the coordinate ring of the $SU(2)$-characters of $\pi_1(H)$, \cite{B}, \cite {PS}. The isomorphism is given by sending a skein given by a disjoint union of simple closed curves $c_i$ to the function that sends a representation $\rho$ to 
\begin{equation}
\prod_i -tr(\rho(c_i)).
\end{equation}
Weyl orthogonality implies that via this isomorphism at $t=-1$ the Yang-Mills measure on a handlebody of genus $g$ is equal to integration on $SU(2)^g$ against Haar measure.

\subsection{Quantum $3$-manifold invariants}

In this section we discuss the evaluation of some quantum invariants of $3$-manifolds.

Suppose that a closed, oriented $3$-manifold $M$ is obtained by surgery on a framed link $L$. Let  $b_+$  (respectively $b_-$) denote the number of positive (negative) eigenvalues of the linking matrix of $L$. The Witten-Reshetikhin-Turaev (WRT) invariant of $M$ at level $r$, where $r>3$ is a positive integer  can be obtained as the following expression:

\begin{equation}
Z_r(M) = <L(\omega)><U_+(\omega)>^{-b_+}<U_-(\omega)>^{-b_-}.
\end{equation}
In this formula $L(\omega)$ denotes decorating each component of a given link with 
$\omega$, the symbol $U_+$ (respectively $U_-$) denotes an unknot with framing $+1$ (respectively $-1$), and $<L>$ denotes taking the value of the Kauffman bracket of a link $L$ in a $3$-sphere evaluated at $4r$-th root of unity.  This formula is taken from Lickorish \cite{Li}. Note that  in this normalization, $Z_r(S^3)=1$ for all $r$.

The Turaev-Viro invariant, $TV(M)$, was originally defined combinatorially in \cite{TV} using triangulations of $3$-manifolds and quantum $6j$-symbols.
Roberts \cite{Rob} proved that the square of the norm of the WRT-invariant of any $3$-manifold $M$ is equal to $TV(M)$. Roberts showed that the TV invariant is equal to the chain-mail invariant derived from a Heegaard diagram of $M$. He used a different normalization than Lickorish. In Robert's normalization $Z_r(S^3)=\frac{1}{X}$, and this agrees with Witten's.

Consider the presentation of the $3$-manifold $M$ by a Heegaard diagram on a standard handlebody $H$ in $S^3$.  Let $\Omega=\frac{1}{X}\omega$. The chain-mail link consists of the attaching curves decorated with $\Omega$ along with the push-offs of a complete set of meridians of $H$ which link the attaching curves and are the $0$-framed unlink with $g$ components, and are also decorated with 
$\Omega$.
The chain-mail invariant is equal to the value of the Kauffman bracket of the chain-mail link multiplied by $\frac{1}{X^2}$.  
Note that
the link consisting of meridians and attaching curves obtained from a Heegaard diagrams yields a surgery presentation for  $M\sharp\overline{M}$.

While the Yang-Mills measure is defined when the complex parameter $t$ has absolute value different than $1$,  the WRT and TV invariants of $3$-manifolds are not defined away from roots of unity. Although the terms of the state sum originally used to define the Turaev-Viro invariant of a manifold $M$ are all well defined when the value of the parameter is equal to $-1$, the sum is infinite and does not converge.

\subsection{Total torsion versus quantum invariants}

By (\ref{globalformula}) the total Dubois torsion is equal to
\begin{equation}\label{glob}
lim_{\lambda\rightarrow \infty} 2^{\frac{3k-3}{2}} \sqrt{\frac{4\pi}{\lambda}}\int_{SU(2)^k}
||\vv||\prod_{i=1}^{k-1}
\frac{1}{2\pi^2}\sum_{c=0}^{\infty}(-1)^{c}(c+1)e^{-c(c+2)/\lambda}s_c(-tr (r_i))\nu^k.
\end{equation}

Recall that integration on $SU(2)^k$ against Haar measure agrees with the Yang-Mills measure on the Kauffman bracket skein module at $t=-1$ of a handlebody of genus $k$. 
Note also that the induced  Riemannan measure on $SU(2)$ as $S^3\subset {\mathbb{R}}^4$ is equal to $\pi^2$ times the Haar measure.
Unfortunately, it is difficult to recognize the integrand in (\ref{glob}) as a skein. Looking at a portion of the formula 
\begin{equation}\label{firstbs}
\sum_{c=0}^{\infty}(-1)^{c}(c+1)s_c(-tr (r_i))= 
\lim_{n\rightarrow\infty}\sum_{c=0}^{n}(-1)^{c}(c+1)s_c(-tr (r_i)),
\end{equation}
and comparing it with (\ref{omega}), note that on the right we have a skein in a handlebody obtained by decorating the attaching curve for a handle corresponding to the $i$-th relator with an analogue of $\omega\in K_{r,f}(S^1\times D^2)$. Thus we can view the left-hand side of (\ref{firstbs}) as {\em coloring} the attaching curves corresponding to  relators $r_i$ with $\omega_{\infty}$. Recall now the method of computing the Yang-Mills measure in $K_{r,f}$ via the Kauffman bracket of a surgery diagram in the 3-sphere. The analogous equation in $K_{r,f}$ for the value of the Yang-Mills measure of the skein 
$\prod_{i=1}^{k-1}\sum_{c=0}^{r-1}(-1)^{c}(c+1)s_c(-tr (r_i))$ is obtained by computing the Kauffman bracket of a chain mail link. The chain mail link is a surgery diagram for the double of the knot complement.

If we tried naively to extend the WRT-invariant from $4r$-th roots of unity to a value at $-1$, 
we could begin by trying to compute the Kauffman bracket of a chain mail link, where instead of decorating the appropriate curves with $\omega$ we used  $\omega_{\infty}$. Unfortunately this infinite sum does not converge. The second attempt would be to regularize using a bump function.   Thus we could interpret the limit
\begin{equation}
lim_{\lambda\rightarrow \infty}  \sqrt{\frac{4\pi}{\lambda}}\int_{SU(2)^k}
\prod_{i=1}^{k-1}
\frac{1}{2\pi^2}\sum_{c=0}^{\infty}(-1)^{c}(c+1)e^{-c(c+2)/\lambda}s_c(-tr (r_i))\nu^k.
\end{equation}
as the regularized WRT-invariant of the double of the knot complement at level $-1$. Alas, this quantity is not a manifold invariant. In order to obtain an invariant we need to multiply the integrand by the factor $||\vv||$. In this way we get an invariant of the double of the knot complement, which is equal to the total Dubois torsion of the knot. In the next section we explain the appearance of $||\vv||$. 

The limit of the WRT-invariants of the double of the knot complement, as the level $r$ tends to infinity, does not exist  due to oscilation. However, the deliberations in this paper lead us to the following conjecture.

\begin{conj}\label{conjecture}
Let $K$ be a knot in $S^3$ whose complement has sufficiently regular $SU(2)$-representation variety. Denote by $T^2$ the boundary torus of the knot complement and by $M$ the $3$-manifold which is the double of the knot complement.
Let $f\in K_{-1}(T^2)$ be a peripheral skein, and let $Z_r(M)$ denote the Witten-Reshetikhin-Turaev invariant of $M$ at level $r$. Then
\begin{equation}
\lim\sup _{r\rightarrow\infty}\frac{1}{r}|Z_r(M,f)|
\end{equation}
defines a seminorm on $K_{-1}(T)$  whose radical is the $A$-ideal of the knot.

\end{conj}

\section{Dubois torsion and the square root of Reidemeister torsion of the double}

In this section we prove that Dubois' torsion of a knot complement is a geometrically motivated evaluation of the square root of the Reidemeister torsion of its double. When the character variety of a knot is nice, the total Dubois' torsion is proportional to the conjectured formula for the leading asymptotics of the Witten-Reshetikhin-Turaev invariant of the double of the knot complement.

 Suppose that  $\rho:\pi_1(M)\rightarrow SU(2)$ is irreducible, and $h$ is a basis
for $H^1(M;ad\rho)$. Choose $h^*$ to be a basis for $H^2(M;ad\rho)$ that
is dual to $h$ with respect to the cup pairing
\begin{equation} \cup:H^1(M;ad\rho)\otimes H^2(M;ad\rho)\rightarrow H^3(M;\mathbb{R}).\end{equation} 
Since $\rho$ is irreducible, 
\begin{equation}
H^0(M;ad\rho)\simeq H^3(M;ad\rho)=\{\vec{0}\}.
\end{equation}
Thus
\begin{equation}
H^*(M;ad\rho)\simeq H^1(M;ad\rho)\oplus H^2(M;ad\rho).
\end{equation}
The {\em square root of the Reidemeister torsion of $M$} evaluated on $\{h,h^*\}$, denoted by
$\sqrt{\tau}(M;h)$, is the square root of the Reidemeister torsion of the chain complex
$C^*(M;ad\rho)$ with respect to the choice of basis $\{h,h^*\}$.   As recognized
by Rozansky \cite{Roz}, $\sqrt{\tau}(M;h)$ is a naturally defined volume on the part of the character variety of $M$ coming from regular irreducible representations.

A {\em volumed vector space} is a vector space $V$ along with a choice
of a nonzero element of the top dimensional exterior power of the vector space,
called the {\em volume}.  The top dimensional exterior power of the direct sum
of two vector spaces is canonically isomorphic to the tensor product of the top dimensional exterior powers of the two vector spaces.

Suppose that $A$, $B$ and $C$ are  chain complexes defined over a field and there is a short
exact sequence
\begin{equation} \begin{CD} 0@>>>A@>\alpha>>B@>\beta>>C@>>>0\end{CD}.\end{equation}
Since the chain
groups are vector spaces, for each $i$ the sequence splits, so that
$B_i\cong A_i\oplus C_i$. We say the sequence is {\em volume exact} if the volume on $B_i$ is the tensor product of the volumes on $A_i$ and $C_i$. There is an approach to Reidemeister torsion due to Milnor \cite{M} based on volumed vector spaces. 
If $\{v_1,\dots, v_2\}$ is the preferred bases for $V$ then 
$v_1\wedge v_2\wedge\dots\wedge v_k$ is the preferred volume.

From the long exact sequence in homology,
we can split the homology groups of $A$, $B$ and $C$.  Let $E=\Delta(J)$,  where
$\Delta$ is the connecting homomorphism and $J$ is a complement to its kernel. Let $F$ be a complement to $E$, and $G=\alpha(F)$. Finally let $H$ be a complement of $G$ and $I=\beta(H)$. With these choices  $H_*(A)=E\oplus F$,  $H_*(B)=G\oplus H$  and $H_*(C)=I\oplus J$.

Choose bases $\mathbf{f}$ for $F$, $\mathbf{h}$ for $H$ and $\mathbf{j}$ for $J$.
Let $\overline{\mathbf{f}}$, $\overline{\mathbf{h}}$, and $\overline{\mathbf{j}}$, denote their images under $\alpha$, $\beta$ and $\Delta$ (on the chain level).  Let $\tau_A(\overline{\mathbf{j}},\mathbf{f})$, $\tau_B(\overline{\mathbf{f}},\mathbf{h})$,$\tau_C(\overline{\mathbf{h}},\mathbf{j})$ denote the Reidemeister torsion of each complex with respect to the indicated choice of bases. The following proposition 
is a consequence of  Theorem 3.2 in \cite{M}, or Corollary 1.2 in \cite{Fre}.

\begin{prop} \label{one}$\tau_B(\overline{\mathbf{f}},\mathbf{h})=\tau_A(\overline{\mathbf{j}},\mathbf{f})\tau_C(\overline{\mathbf{h}},\mathbf{j}) $\end{prop}
 \qed

In order to analyze the cohomology of the double we will use the Mayer-Vietoris sequence coming from the fact that the double is the union of two copies  of the knot complement. 
Suppose that $M=A\cup B$ is excisive with inclusion map $i:A\cap B\rightarrow M$, so that there are two Mayer-Vietoris sequences, one with coefficients $G$,
\begin{equation}
\begin{CD}H^{i-1}(A\cap B;G)@>\Delta>>H^i(M;G)@>>>H^i(A;G)\oplus H^i(B:G) @>>>H^i(A\cap B;G)\end{CD},
\end{equation}
and the other with coefficients $G'$,
\begin{equation}
\begin{CD}H^{i-1}(A\cap B;G')@>\Delta'>>H^i(M;G')@>>>H^i(A;G')\oplus H^i(B:G') @>>>H^i(A\cap B;G')\end{CD}.
\end{equation}
Also suppose that there is a symmetric
pairing $G\otimes G\rightarrow G'$ that allows us to define a cup pairing,
\begin{equation} 
\cup:H^i(M;G)\otimes H^j(M;G)\rightarrow H^{i+j}(M;G').
\end{equation}
For any $\alpha \in H^i(A\cap B,G)$ and $\beta\in H^j(M;G)$ we have
\begin{equation} \Delta(\alpha)\cup \beta=\Delta'(\alpha \cup i^*\beta).\end{equation}
This formula can be found on page 252 of \cite{Sp} except for the use of local coefficients.
However,  local coefficients can always be reduced to standard coefficients, by passing to a cover and taking the equivariant part, so the result holds in our setting.

Suppose now that $S^3-K$ is a knot complement. For this part of the paper we mean that we have removed an open regular neighborhood of $K$ from $S^3$, so that $S^3-K$ is a compact $3$-manifold with boundary a torus $T^2$.  Denote by $M$  the double
of $S^3$, so that $M=X_1\cup X_2$, where $X_1$ and $X_2$ are two copies of $S^3-K$, with the orientation of $X_2$  the opposite from the orientation on $S^3-K$, and $X_1$ and $X_2$ identified along their boundaries via the identity map.  There is a regular CW decomposition of $S^3-K$ with one vertex $v$ on the boundary, which can be doubled to get a CW-decomposition of $M$.

In order to relate the torsion of the knot complement to the square root of the torsion of the double we need to choose the bases for the cohomology groups for all the spaces in the Mayer-Vietoris sequence for $M=X_1\cup X_2$.

Let
$\rho$ be a regular $SU(2)$ representation of $S^3-K$. Let $\mu$ denote the meridian of $S^3-K$ and
$\rho(\mu)=\cos\phi +sin\phi \vP$.  Let $g=exp(\theta \vP)$, where $0\leq\theta\leq\pi$.  
Let
\begin{equation} \rho_g:\pi_1(M)\rightarrow SU(2)\end{equation}
be the representation that restricts to $\rho$ on $\pi_1(X_2)$ and to $g\rho g^{-1}$ on $\pi_1(X_1)$. The construction of a path of representations depending on the variable $\theta$ like this is called {\em bending}. The tangent vector to such a path can be thought of as a tangent vector to a {\em deformation}.
 Note that $C^*(S^3-K;ad(g\rho g^{-1}))$ is canonically isomorphic to $C^*(S^3-K;ad\rho)$, under the isomorphism $c^*\otimes \vv\mapsto c^*\otimes g.\vv$. By dint of this isomorphism we identify $C^*(X_1;ad(g\rho g^{-1}))$ and $C^*(X_2;ad\rho)$ with $C^*(S^3-K;ad\rho)$.

Since $\rho$ restricted to the boundary is not central, $H^*(T^2;ad\rho)=H^*(T^2;{\mathbb {R}})\otimes {\mathbb{R}}\vP$. A volume basis for  $H^*(T^2;ad\rho)$ can be chosen as follows. Denote by $v^*$ a cochain in $C^0(T^2;{\mathbb{R}})$ dual to the vertex $v$, and let $t\in C^2(T^2;{\mathbb{R}})$ be a cochain with $t[T^2]=1$. Finally, let $\xi$ and $\eta$ be the cocycles in $ C^1(T^2;{\mathbb{R}})$ such that $\xi\cup\eta=t$. A volume basis consists of $v^*\otimes \vP$, $\xi\otimes \vP$, $\eta\otimes\vP$ and $t\otimes\vP$.

Since $\rho$ is regular $H^1(S^3-K;ad\rho)$ is one-dimensional and is spanned by the tangent vector $h$ to a path of representations. The restriction map $H^2(S^3-K;ad\rho)\rightarrow H^2(T^2;ad\rho)$ is an isomorphism. Finally $H^0(S^3-K;ad\rho)=0$ as $\rho$ is irreducible.

From this we conclude that $H^1(M;ad\rho_g)$ is two-dimensional. One dimension comes from the image of $H^0(T^2;ad\rho)$ under the connecting homomorphism $\Delta$. The second is accounted for by the fact that any complement of the image of $\Delta$ is mapped injectively into the diagonal of the direct sum 
\begin{equation} H^1(X_1; ad (g \rho g^{-1}))\oplus H^1(X_2;ad\rho)= H^1(S^3-K;ad\rho)\oplus H^1(S^3-K;ad\rho).
\end{equation}
The second cohomology of $M$ is also two-dimensional where analogously one dimension is accounted for by the image of the connecting homomorphism and any complement of that image gets mapped injectively into the diagonal of
\begin{equation} H^2(X_1;ad(g\rho g^{-1}))\oplus H^2(X_2;ad\rho)= H^2(S^3-K;ad\rho)\oplus H^2(S^3-K;ad\rho).
\end{equation}

There is a nondegenerate pairing coming from the cup product on $M$ and the dot product on the coefficients,
\begin{equation} \cup :H^1(M;ad\rho_g)\otimes H^2(M;ad\rho_g)\rightarrow H^3(M;\mathbb{R}).\end{equation}
It should be noted that if $m\in H^3(M;\mathbb{R})$ is a cohomology class with $m([M])=1$ then
$\Delta'(t)=m$ where $\Delta'$ is the connecting homomorphism from the  Mayer-Vietoris sequence for $M=X_1\cup X_2$ with real coefficients.

If $t'\in H^2(S^3-K;ad\rho)$ is chosen so that
restriction of $t'$ to $H^2(T^2;ad\rho)$ is $t$ and $h$ is a basis for $H^1(S^3-K;ad\rho)$ then the Dubois torsion  $\tau(S^3-K;h)$ is the Reidemeister torsion of $C^*(S^3-K;ad\rho)$ with respect to the bases $h$ for $H^1(S^3-K;ad\rho)$ and $t'$ for $H^2(S^3-K;ad\rho)$.  

\begin{theorem} \label{comparison} Let $M$ be the double of $S^3-K$ and let $\rho_g:\pi_1(M)\rightarrow SU(2)$ be obtained from bending  $\rho:\pi_1(S^3-K)\rightarrow SU(2)$ as above. Assume that $\rho$ and $\rho_g$ are regular.  If we choose $h$ as the basis for $H^1(S^3-K;ad\rho)$ and we use the basis $\{\Delta(v^*\otimes \vP),\tilde{h}\}$ for $H^1(M;ad\rho)$ where $\tilde{h}$ restricts to $h$, then
\begin{equation} \sqrt{\tau}(M;\{\Delta(v^*\otimes \vP),\tilde{h}\})=\tau(S^3-K;h).\end{equation} \end{theorem}

\proof The proof is an application of Proposition \ref{one} followed  by some interpretation.
We apply the proposition to the short exact sequence inducing the Mayer-Vietoris sequence
\begin{equation}\label{shmv}
0\rightarrow C^*(M;ad\rho_g)\rightarrow C^*(X_1;ad(g\rho g^{-1}))\oplus C^*(X_2;ad\rho)\rightarrow C^*(T^2;ad\rho)\rightarrow0.
\end{equation}
These are the chain complexes associated to a regular $CW$-decomposition of $S^3-K$ with a single vertex $v$ on the boundary, which we double to get a $CW$-decomposition of $M$. We use the  innerproduct on $C^*(M;ad\rho)$,  coming from the CW structure and 
 the dot product on $su(2)=\mathbb{R}^3$ .
The sequence (\ref{shmv}) is volume exact.  

Denote by $\tilde{h}$  a cocycle in $C^1(M;ad\rho)$ that restricts to the basis $h$ chosen for $H^1(S^3-K;ad\rho)$.
  We use the basis $\{\Delta(v^*\otimes \vP), \tilde{h}\}$ for $H^1(M;ad\rho_g)$. 
The pair  $\{(h,h),(h/2,-h/2)\}$  forms a basis for 
\begin{equation}
H^1(X_1;ad(g\rho g^{-1}))\oplus H^1(X_2;ad\rho)=H^1(S^3-K;ad\rho)\oplus H^1(S^3-K;ad\rho)
\end{equation}
Notice that the vector $(h,h)$  is the image of $\tilde{h}$ under the map in the Mayer-Vietoris sequence.  Denote the image of $(h/2,-h/2)$ under the difference of the restriction maps by $t_1$.
Note that $t_1$
is a tangent vector to the image under restriction of the curve of regular characters $X^i(S^3-K)$ in the pillowcase. 

Let $n$ be a normal vector to $X^i(S^3-K)$ at $[\rho]$, so that $\{t_1,n\}$ is a volume basis for $H^1(T^2;ad\rho)$.
Let $h'\in C^2(S^3-K;ad\rho)$ be a cocycle that restricts to $t\otimes \vP$ in $H^2(T^2;ad\rho)$.
We use $\{(h',h'),(h'/2,-h'/2)\}$ as a basis for 
\begin{equation}
H^2(X_1;ad(g\rho g^{-1}))\oplus H^2(X_2;ad\rho)=H^2(S^3-K;ad\rho)\oplus H^2(S^3-K;ad\rho).
\end{equation}  
Let
$h_2\in C^2(M;ad\rho_g)$ be a cocycle that pulls back to $(h',h')$.  Finally use $\{\Delta(n),h_2\}$
as the basis for $H^2(M;ad\rho_g)$.

We have chosen our bases as in the hypothesis of Proposition \ref{one}. Thus
\begin{equation}\tau(M;\{\Delta(v^*\otimes \vP), \tilde{h},\Delta(n),h_2\})=\tau(T^2,\{v^*\otimes \vP,t_1,n,t\})\times\end{equation}
\[\tau\left(C^*(X_1;ad\rho)\oplus C^*(X_2;ad\rho);\{(h,h),(h/2,-h/2),(h',h'),(h'/2,-h/2)\}\right).\]

We chose the basis for $T^2$ to be a volume basis, so
\begin{equation} \tau(T^2,\{v^*\otimes \vP,t_1,n,t\})=1.\end{equation}

The basis for $H^*(M;ad\rho_g)$ was chosen so that it computes \begin{equation}\sqrt{\tau}(M;\{\Delta(v^*\otimes \vP), \tilde{h}\}).\end{equation}  To see this we need to compute the three cup pairings \begin{equation}\Delta^0(v^*\otimes \vP)\cup \Delta^1(n), \Delta^0(v^*\otimes \vP)\cup h_2,\  \mathrm{and} \ \tilde{h} \cup \Delta^1(n).\end{equation}  We do this as follows:
\begin{enumerate}
\item 
\begin{equation}\Delta^0(v^*\otimes \vP)\cup \Delta^1(n)=\Delta^1(v^*\otimes \vP\cup i^*\Delta^1(n))=0\end{equation}
as $i^*\Delta^1(n)=0$,
\item
\begin{equation}\Delta^0(v^*\otimes \vP)\cup h_2=\Delta^2(v^*\otimes \vP\cup i^*h_2)=1\end{equation}
as $i^*h_2=t\otimes \vP$, and
\item
\begin{equation} \tilde{h} \cup \Delta^1(n)=\Delta^1(n)\cup \tilde{h}=\Delta^2(n\cup i^*\tilde{h})=-1\end{equation}
as $i^*\tilde{h}=t_1$.

\end{enumerate}

We don't care what the last pairing is because we can make a volume preserving change of basis to get a dual basis.

Finally, we interpret
\begin{equation}\tau(C^*(X_1;ad\rho)\oplus C^*(X_2;ad\rho);\{(h,h),(h/2,-h/2),(h',h'),(h'/2,-h'/2)\}).\end{equation}
Change bases by elementary row operations that do not change volume to
get
\begin{equation}\tau(C^*(X_1;ad\rho)\oplus C^*(X_2;ad\rho);\{(h,h),(h/2,-h/2),(h',h'),(h'/2,-h'/2)\})=\end{equation}
\begin{equation}\tau(C^*(X_1;ad\rho)\oplus C^*(X_2;ad\rho);\{(h,0),(0,-3h/2),(h',0),(0,-3h'/2)\}).\end{equation}
The torsion now splits as the product of torsions, the first computes the Dubois torsion
of $S^3-K$ at $\rho$.  The second also computes Dubois torsion as  $-3/2$ appears once in the bottom and once in the top so it cancels.  Hence we get
the desired result. \qed

To understand this geometrically, consider
a path of representations covering an arc $Y(S^3-K)$ of regular representations in the character variety of the knot complement $S^3-K$.  By bending as above we can parametrize a circle bundle over the arc $Y(S^3-K)$ which we call $Y(M)$. There is a projection that just restricts the representation to the second copy of $S^3-K$, call it
\begin{equation}\sigma:Y(M)\rightarrow Y(S^3-K).\end{equation}
The kernel of the derivative of $\sigma$ is the span of $\Delta^0(v^*\times \vP)$.  The partial derivative of the parametrization in the first variable is $\Delta^0(v^*\times \vP)$.

\begin{prop} The measure on $Y(S^3-K)$ coming from Dubois torsion is $\frac{1}{\pi}$ times the pushforward measure under $\sigma$ of the measure on $Y(M)$ coming from the square root of Reidemeister torsion. \end{prop}

\proof The volume of the fiber is
$\pi||\Delta^0(v^*\times \vP)||$ which is $\pi||\vv||$.  If $S\subset Y(S^3-K)$ is measurable,
then the measure of $\sigma^{-1}(S)$ is $\pi$-times the measure of $S$. \qed

The factor $||\vv||$ that appears in the formula for Dubois torsion can be attributed to partial integration over the fiber of the map $\sigma$.

Assume now that $K\subset S^3$ is a knot so that the regular representations are dense in the nonabelian representations, and the map $X(S^3-K)\rightarrow X(T^2)$ is an embedding away from finitely many points. The character variety of $M$ has a component coming from  bent representations lying over each component coming from regular representations of $\pi_1(S^3-K)$.  

\begin{prop} Let $Y^i(M)$ be the part of the character variety of $M$ that comes from bending regular representations of $\pi_1(S^3-K)$.  Let
$f:X^i(S^3-K)\rightarrow \mathbb{R}$ be integrable, then
\begin{equation} \int_{Y^i(M)} f \sqrt{\tau}(M)=\pi \int_{X^i(S^3-K)}f \ d\tau\end{equation} where we lift $f$ to $Y^i(M)$ via the projection map.  \end{prop}\qed

Recall the conjectured formula for the leading asymptotics of the WRT-invariant,

\begin{equation}
Z_r^{SU(2)}(M)\sim_{r\rightarrow \infty} e^{-3\pi \bf{i}(1+b^1(M))/4}\times 
\int_{[A]}e^{2\pi \bf{i}CS(A)} r^{(h^1_A-h^0_A)/2} e^{-2\pi\bf{i}(I_A/4+h^0_A/8)}\tau_M(A)^{1/2}.
\end{equation}

We consider the right hand side of the formula at characters of representations of $M$
that come from bending regular representations of $\pi_1(S^3-K)$.
With this restriction $h^1=2$, $h^0=0$, $b^1(M)=1$. Since $M$ is the double of a knot complement the Chern-Simons invariant is an integer and the spectral flow is zero \cite{KK,KL}, so along the regular representations the absolute value of the  right hand side reduces to 

\begin{equation} 
r\int_{[A]\in X^i(M)}  \tau_M(A)^{1/2}.
\end{equation}

If to naively assume that there were no contributions to the leading order asymptotics coming from the singular points of the character variety, then the total Dubois torsion of $S^3-K$ is  proportional to leading order asymptotics of the WRT-invariants of the double.

Comparing this now to  Conjecture \ref{conjecture}, we are working with the conjectural asymptotic formula extended to the WRT-invariant for a knot that lies in the boundary of the knot complement inside of the double of the knot complement. 
In the case of knots with very well behaved character varieties our conjecture is a consequence of the conjectured asymptotics of the WRT-invariant.

\end{document}